%% file: main.tex
\documentclass[a4paper, 11pt, reqno]{amsart}
\usepackage{subfiles} % foldered structure
\usepackage{caption} % caption width
\usepackage{hyperref}
\usepackage[export]{adjustbox} % frame around figure

\hypersetup{colorlinks=true,urlcolor=blue,citecolor=red}
\newtheorem{theorem}{Theorem}[section]

\newtheorem{definition}[theorem]{Definition}

\newtheorem{remark}[theorem]{Remark}

\begin{document}

\title{A Lighthouse Illumination Problem}
\date{March 3, 2019}
\author{Erhan Tezcan}

\begin{abstract}
This paper discusses a problem that consists of $n$ ``lighthouses'' which are circles with radius 1, placed around a common center, equidistant at $n$ units away from the placement center. Consecutive lighthouses are separated by the same angle: $360^\circ/n$ which we denote as $\alpha$. Each lighthouse ``illuminates'' facing towards the placement center with the same angle $\alpha$, also called ``Illumination Angle'' in this case. As for the light source itself, there are two variations: a single point light source at the center of each lighthouse and point light sources on the arc seen by the illumination angle for each lighthouse. The problem: what is the total dark (not illuminated) area for a given number of lighthouses, and as the number of lighthouses approach infinity? We show that by definition of the problem, neighbor lighthouses do not overlap or be tangent to each other. We propose a solution for the center point light source case, and discuss several small cases of $n$ for the arc light source case. 
\end{abstract}

\maketitle

\subfile{Sections/1-Introduction}

\subfile{Sections/2-Definition}

\subfile{Sections/2.1-Boundaries}

\subfile{Sections/2.2-Center-Light-Source}

\subfile{Sections/2.3-Arc-Light-Source}

\subfile{Sections/3-Conclusion}

\end{document}

%% file: Sections/1-Introduction.tex
%% BEGIN %%
\section{Introduction}
This problem originated to the author during a 17-hour long bus trip from Warsaw to Tallinn, in line with author's inability to sleep in a bus. At night, with only visual input being the bus' indoor ceiling lights that are \textbf{circular}, the author tried to pass time by thinking of a random problem originating from aforementioned light sources. This paper discusses that problem, which the author initially called ``Lighthouse Problem'' but to differentiate it from another problem of the same name that takes place in \cite{lighthouse-problem}, we would like to call our problem ``A Lighthouse Illumination Problem''. \par
In the next section we will be defining the problem, afterwards, we will discuss the boundaries and then delve into the problem itself. Finally, we will be recapitulating our questions and question marks. We created all our figures using a software called ``GeoGebra'' \cite{geogebra}.
%% END %%

%% file: Sections/2-Definition.tex
%% BEGIN %%
\section{The Problem}
\label{SECT: The Problem}

\begin{figure}
    \centering
    \begin{minipage}{0.45\textwidth}
        \centering
        \captionsetup{width=0.9\textwidth}
        \includegraphics[width=0.9\textwidth, frame]{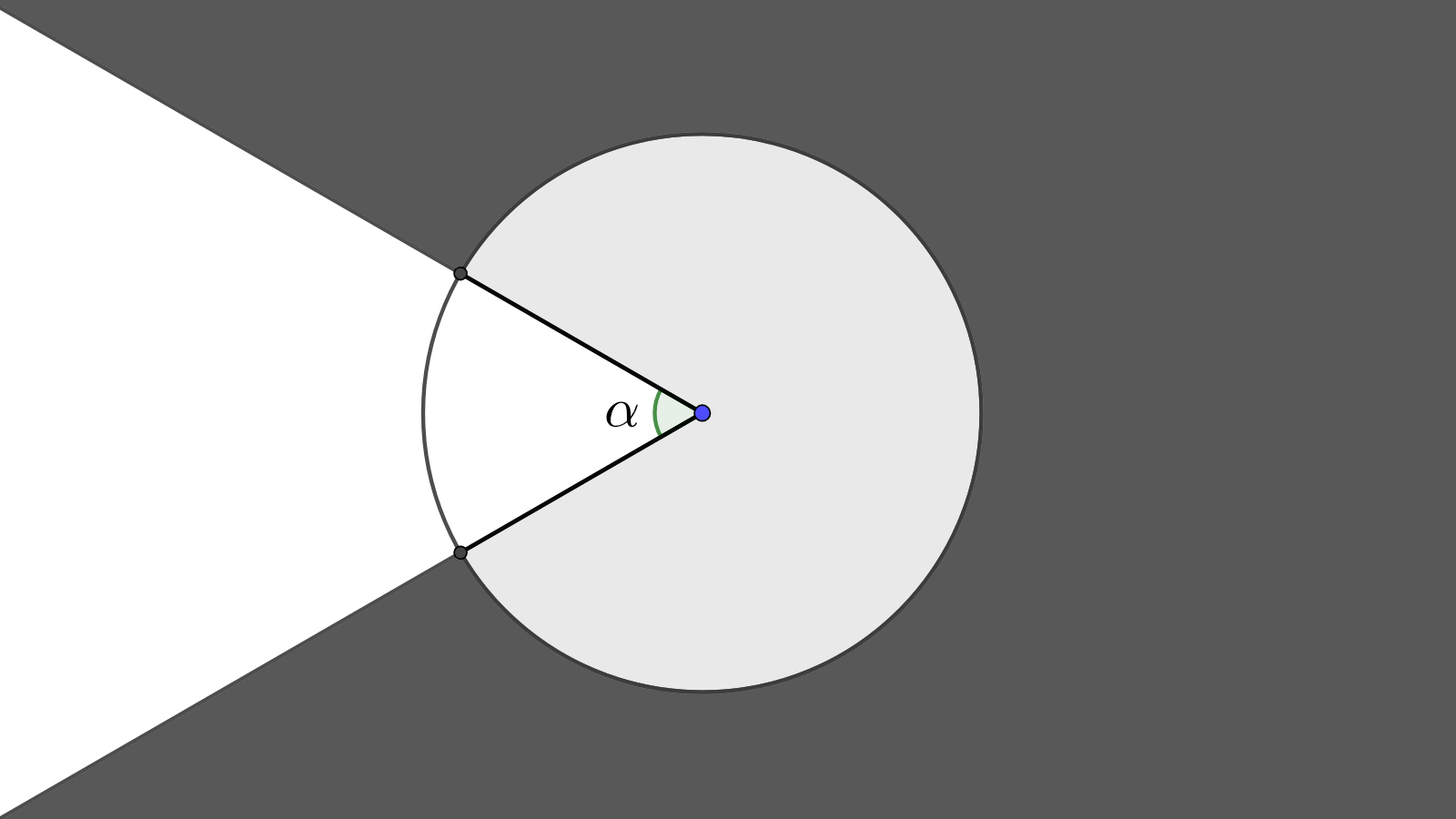}
        \caption{Point Light Source at Center}
        \label{FIG: Center Light Source}
    \end{minipage}
    \hfill
    \begin{minipage}{0.45\textwidth}
        \centering
        \captionsetup{width=0.9\textwidth}
        \includegraphics[width=0.9\textwidth, frame]{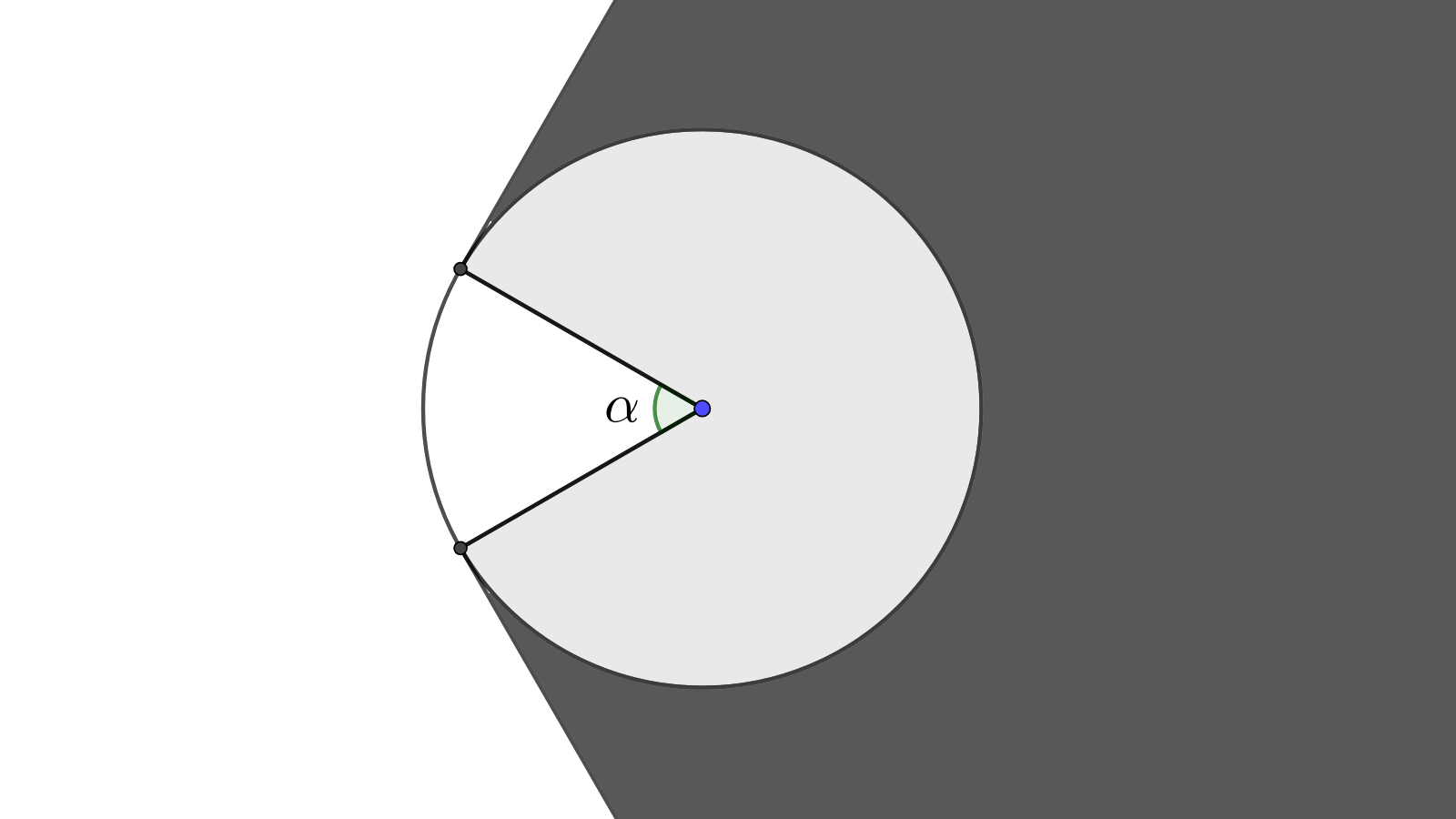}
        \caption{Point Light Source at Arc / Edge}
        \label{FIG: Arc Light Source}
    \end{minipage}
\end{figure}

\begin{definition}[The Lighthouse Illumination Problem]
Suppose we have $n$ circles with radius $1$ in an infinite plane, placed around a common center point (placement center) where the distance between any circle's center and the placement center is $n$. Placing $n$  circles like this divides the $360^\circ$ of the center point into $n$ angles of $360^\circ / n$. We will denote this angle as $\alpha$. Each circle acts as a ``lighthouse'', illuminating towards the center, looking directly at the placement center point with the illumination angle $\alpha$. The light source has two variations in this problem: 
\begin{enumerate}
    \item The center point of each lighthouse is a point light source. (Figure \ref{FIG: Center Light Source})
    \item Every point on the intercepted minor arc of the illumination angle $\alpha$ act as a point light source. (Figure \ref{FIG: Arc Light Source})
\end{enumerate}
For both variations of the light source, what is the total dark area for a given number of lighthouses $n$?
\end{definition}
To further explain figures \ref{FIG: Center Light Source} and \ref{FIG: Arc Light Source}, we shall explain the three shades: white, dark gray and light gray.
\begin{itemize}
    \item White denotes an illuminated area. The circular sector defined by the illumination angle $\alpha$ is also white.
    \item Dark gray denotes a dark area, these are not illuminated by any lighthouse.
    \item Light gray denotes the non-illuminating part of the lighthouse. 
\end{itemize}
One might wonder whether the light gray counts as a dark area or not. We made our calculations with the light gray area excluded, because it is easy to include it if we want to. The light gray area for every lighthouse is given by
\begin{equation}
    \label{EQ: Light Gray Shaded Area Single}
    \pi - \pi \frac{\alpha}{360^\circ} = \pi \left( 1 - \frac{\frac{360^\circ}{n}}{360^\circ}\right) = \pi \left(1 - \frac{1}{n}\right)
\end{equation}
Multiplying \eqref{EQ: Light Gray Shaded Area Single} by $n$ gives the total light gray area
\begin{equation}
    \label{EQ: Light Gray Shaded Area Total}
    n\pi \left(1 - \frac{1}{n}\right) = \pi(n - 1)
\end{equation}
If we want to include the light gray area as darkness we just add $\pi(n - 1)$ to our calculation. \par
By definition, the only possible dark area is the area where no lighthouse can illuminate. The only ``objects'' in the plane are lighthouses, and they very well illuminate in front of them. Their sides are illuminated by the other lighthouses, but behind a lighthouse is not illuminated. The dark area occurs behind a lighthouse. This problem has identical lighthouses, same radius, same distance to placement center, same illumination angle. We can find the total dark area just by looking at the dark area behind one lighthouse, then multiply that with the number of lighthouses $n$. As for notation, we will use $D(n)$ for total dark area, $d(n)$ for the dark area behind a single lighthouse.
\begin{equation}
    \label{EQ: Total Dark Area using d(n)}
    D(n) = n \times d(n)
\end{equation}
\par The dark area behind a lighthouse is defined by two light rays coming from two other lighthouses. We will call the lighthouse we are calculating the dark area behind of as the ``target lighthouse''. The lighthouse bearing the light source of the aforementioned light ray will be called ``source lighthouse''.
%% END %%

%% file: Sections/2.1-Boundaries.tex
%% BEGIN %%
\section{Boundary Cases}
Before we actually delve into the dark area calculations, we want to discuss two cases. For some number of lighthouses $m$ is it possible that
\begin{enumerate}
    \item Lighthouses are tangent to each other?
    \item Lighthouses are overlapping?
\end{enumerate}

\begin{theorem}
    Neighbor lighthouses never overlap or touch for $n \geq 2$.
\end{theorem}
\begin{proof}
We can answer both by looking at the triangle defined by a pair of neighbor lighthouse centers and the placement center. 
\begin{figure}
    \centering
    \begin{minipage}{0.48\textwidth}
        On figure \ref{FIG: Distance Between Lighthouses} we have three points: $L_i$ is a lighthouse, $L_{i+1}$ is the lighthouse placed next to it and $P$ is their common center which they are placed around of. Recall that $\alpha = 360^\circ / n = 2\pi / n$. This figure is true for any pair of neighbor lighthouses for any $n$. 
    \end{minipage}\hfill 
    \begin{minipage}{0.48\textwidth}
        \centering
        \captionsetup{width=0.7\textwidth}
        \includegraphics[width=0.7\textwidth]{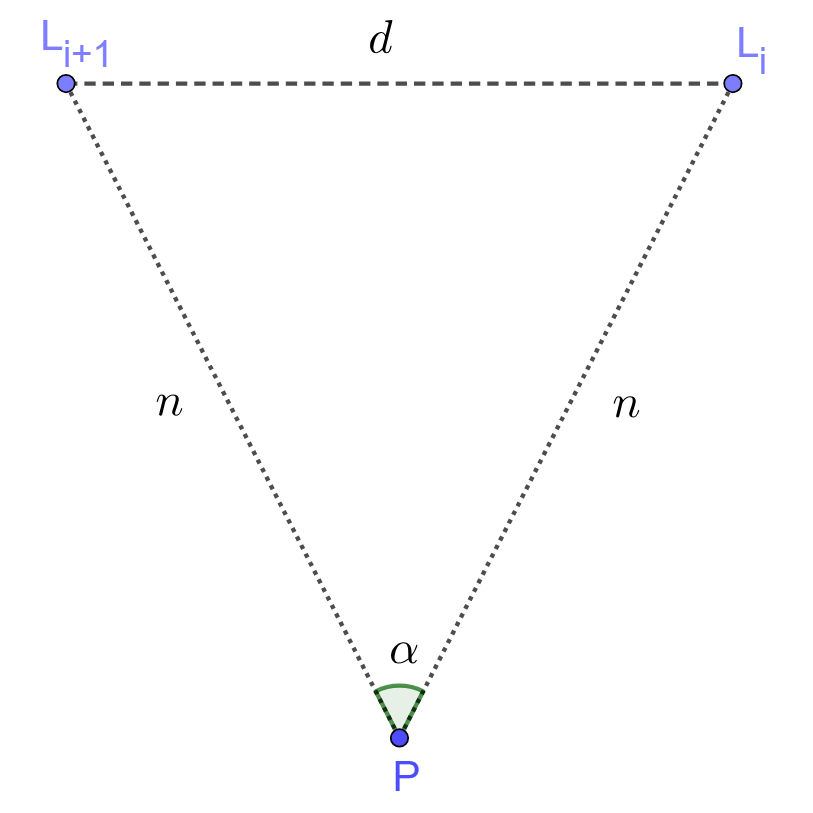}
        \caption{Distance between 2 lighthouses.}
         \label{FIG: Distance Between Lighthouses}
    \end{minipage}
\end{figure}
Each lighthouse has a radius of $1$. Looking at figure \ref{FIG: Distance Between Lighthouses}, if the neighbor lighthouses are tangent, then $d=2$. Similarly, if the neighbor lighthouses are overlapping then $d<2$. We can write $d$ in terms of $n$ using Cosine Theorem. 
\begin{equation*}
    d^2 = n^2 + n^2 - 2n^2\cos(\alpha)
\end{equation*}
This reduces to $d^2 = 2n^2 (1 - \cos(\alpha))$. Recall that $1 - \cos(\phi) = 2\sin^2(\phi)$ so we get $d^2 = 2n^2(2\sin^2(\alpha)) = 4n^2\sin^2(\alpha)$. Taking the root yields
\begin{equation}
    \label{EQ: Cosine Theorem d by n}
    d = 2n \sin(\pi/n)
\end{equation}
Now we can prove the cases by showing that $d > 2$. Note that $d$ only exists when there are 2 or more lighthouses, so $n\geq 2$ and we want $d > 2$ which gives the inequality
\begin{equation*}
    2n \sin(\pi/n) > 2
\end{equation*}
Looking at $n=2$, we have 
\begin{equation*}
    4\sin(90^\circ) = 4 > 2
\end{equation*}
which is correct. From this point on, $2n \sin(\pi/n)$ is monotonically increasing, and since the first value was $4$ the remaining values will never be less than $4$. 
\end{proof}
\begin{remark}
The case of tangent lighthouses $d=2$ can be also studied by referring to Steiner Chain \cite{wolfram-steiner-chain}. Regarding a Steiner Chain:
\begin{equation}
\label{EQ: Steiner Formula}
    \sin\left(\frac{\pi}{n}\right) = \frac{a-b}{a+b}
\end{equation}
where there are $n$ circles packed between a central circle of radius $b$ and an outer concentric circle of radius $a$. Relating this to our problem, we have $b=n-1$ and $a=n+1$. Plugging them in equation \eqref{EQ: Steiner Formula} yields:
\begin{equation*}
    \sin\left(\frac{\pi}{n}\right) = \frac{n+1-n+1}{n+1+n-1} = \frac{2}{2n} = \frac{1}{n}
\end{equation*}
This is the equation we would get by plugging $d=2$ in equation \eqref{EQ: Cosine Theorem d by n}. 
\end{remark}

\subsection{Distance between neighbor lighthouses at infinity.} We can find the distance between neighbor lighthouses at infinity just by looking at 
\begin{equation*}
    \lim_{n\to\infty} 2n\sin(\pi/n) = \lim_{n\to\infty} 2\pi \frac{\sin(\pi/n)}{\pi/n}
\end{equation*}
Substituting $\theta = \pi/n$ gives us
\begin{equation*}
    \lim_{\theta\to0} 2\pi \frac{\sin(\theta)}{\theta} = 2\pi
\end{equation*}
This means that as $n$ approaches infinity the distance between neighbor lighthouses $d$ becomes $2\pi$. Now that we have shown lighthouses will remain unscathed without colliding with one another, we can continue with the dark area calculations.
%% END %%

%% file: Sections/2.2-Center-Light-Source.tex
%% BEGIN %%
\section{Point Light Source at the Center}
Our first variation of the problem has each lighthouses having a single point light source at their center. We will be giving examples for 1, 2, 3, 4, 5 and 6 lighthouses. We will then discuss the dark area for any $n$. 
\begin{figure}
    \centering
    \begin{minipage}{0.33\textwidth}
        \centering
        \captionsetup{width=.9\textwidth}
        \includegraphics[width=0.9\textwidth, frame]{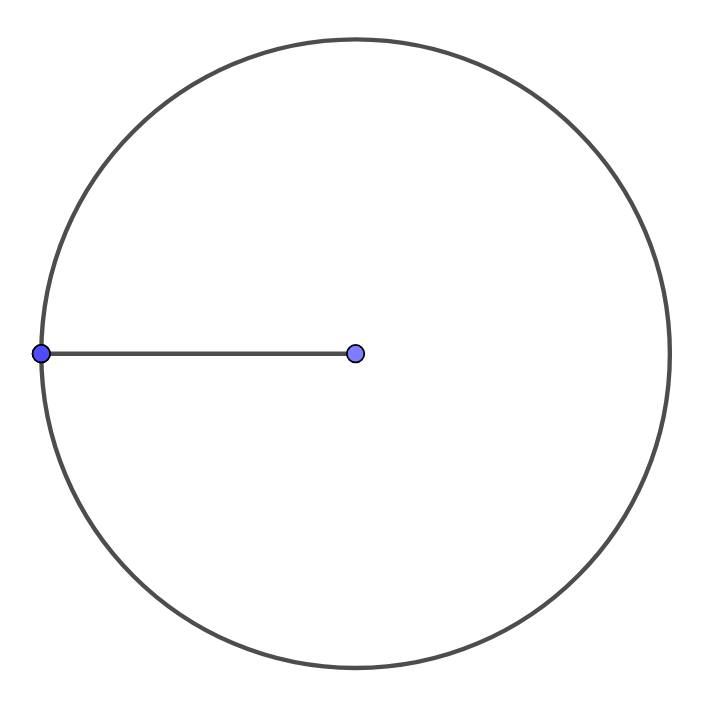}
        \caption{1 Lighthouse.}
         \label{FIG: 1 Lighthouse}
    \end{minipage}\hfill 
    \begin{minipage}{0.67\textwidth}
        \centering
        \captionsetup{width=.9\textwidth}
        \includegraphics[width=0.9\textwidth, frame]{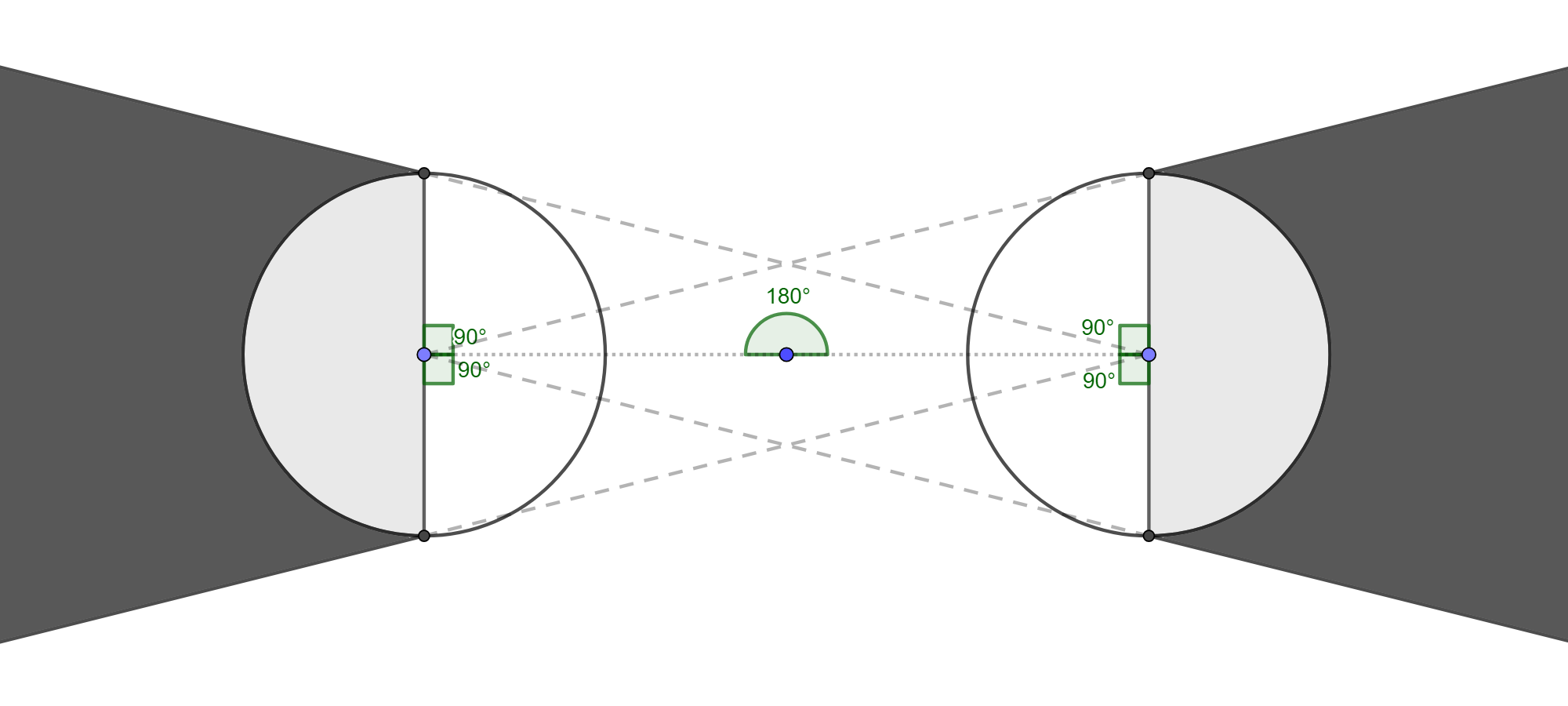}
        \caption{2 Lighthouses, Center Point Light Source.}
        \label{FIG: 2 Lighthouses - Center}
    \end{minipage}
\end{figure}
%%%%%%%%%%%%%%%%%%%%%%%%%%%%%%%%%%%%%%%%%%%%%%%%%%%%%%%%%%%%%%%%%%%%%%
\subsection{1 Lighthouse} 
Our first case is a single lighthouse. Illumination angle $\alpha$ was defined to be $360^\circ/n$. Now that $n=1$ we have $360^\circ$ illumination angle. This basically describes a single lighthouse illuminating in all directions. On figure \ref{FIG: 1 Lighthouse} we can see two points. The point on the left is the placement center, point on the right is the center of the lighthouse. Total dark area is $0$, in other words, $D(1) = 0$.
%%%%%%%%%%%%%%%%%%%%%%%%%%%%%%%%%%%%%%%%%%%%%%%%%%%%%%%%%%%%%%%%%%%%%%
\subsection{2 Lighthouses}
On figure \ref{FIG: 2 Lighthouses - Center} we see the case of $n=2$. Two lighthouses are facing each other and evidently they can not illuminate behind each other. The result is an infinite dark area, $D(2) = \infty$.
%%%%%%%%%%%%%%%%%%%%%%%%%%%%%%%%%%%%%%%%%%%%%%%%%%%%%%%%%%%%%%%%%%%%%%
\subsection{3 Lighthouses}
The first visually appealing case is $n=3$. It is also the first non-zero finite value for the dark area. 
\begin{figure}
    \centering
    \includegraphics[width=0.70\textwidth, frame]{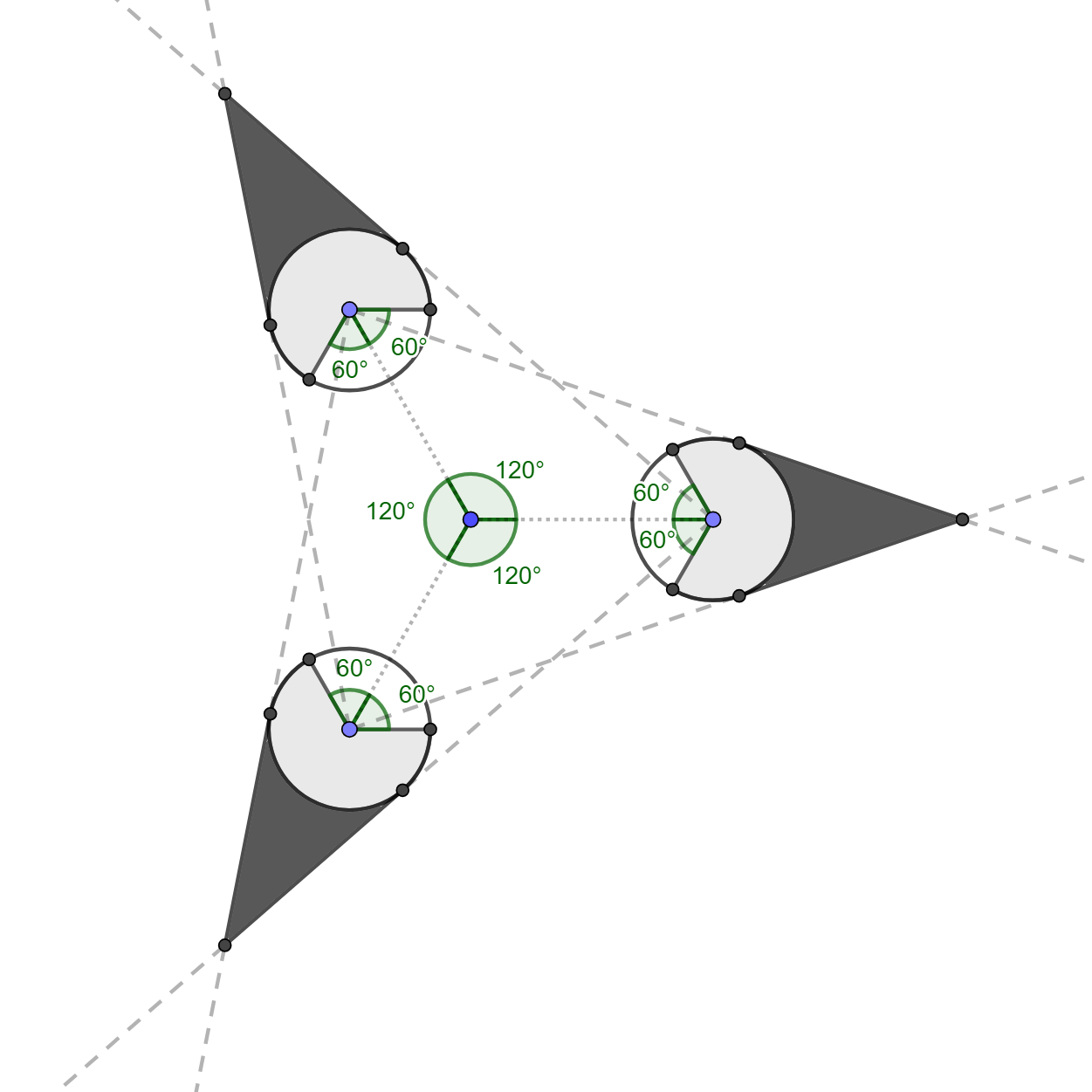}
    \caption{3 Lighthouses, Center Point Light Source.}
    \label{FIG: 3 Lighthouses - Center}
\end{figure}
Our eyes could measure $0$ and $\infty$ but now we will have to do some calculations. We show the $n=3$ case on figure \ref{FIG: 3 Lighthouses - Center}. Recall that in section \ref{SECT: The Problem} we talked about the fact that it is possible to focus on just a single lighthouse and the dark area behind it. That is exactly what we will be doing.
\begin{figure}
    \centering
    \includegraphics[width=0.75\textwidth]{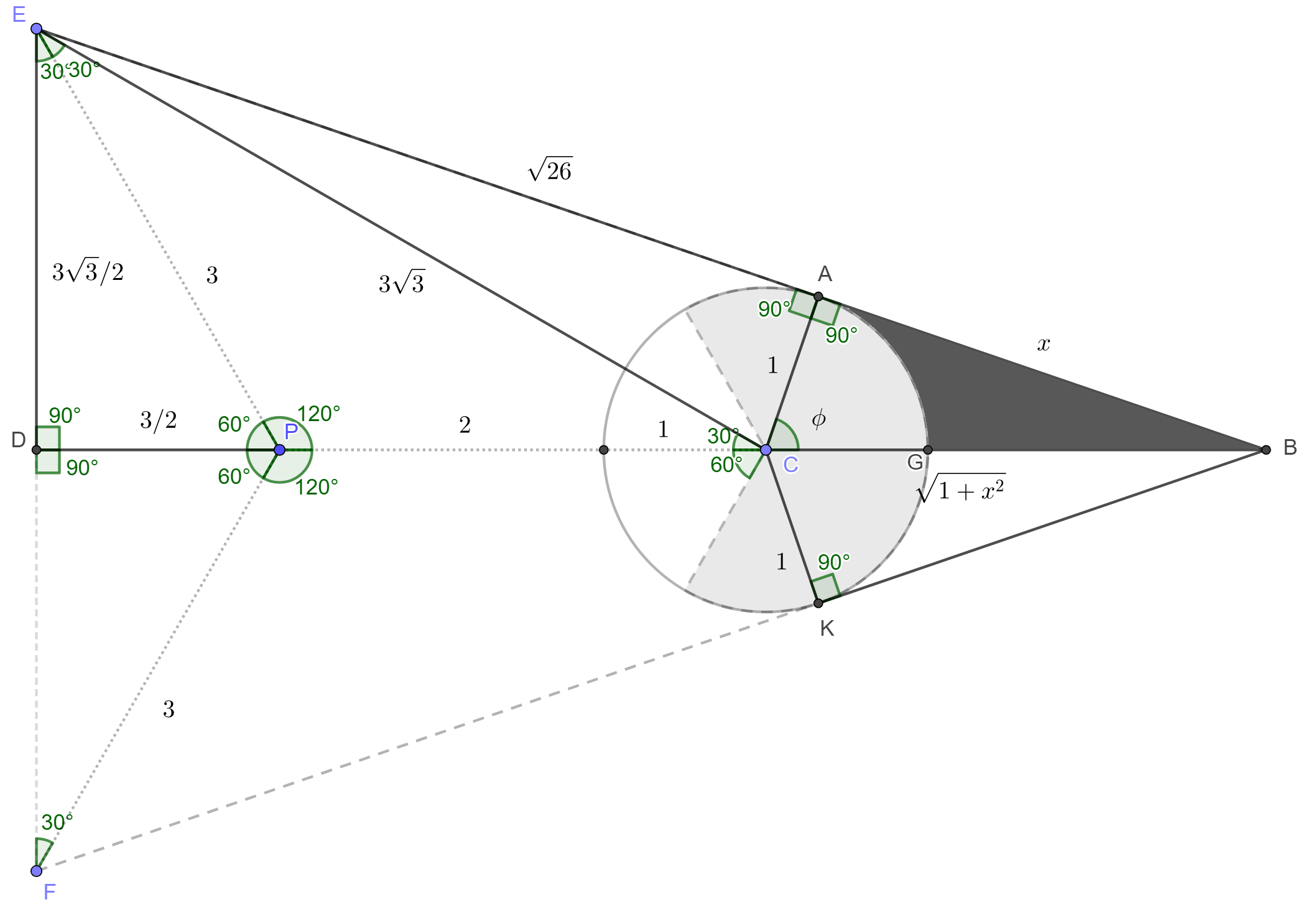} 
    \caption{3 Lighthouses, zoomed in on a target lighthouse.}
    \label{FIG: 3 Lighthouses - Center - Zoom}
\end{figure}
On figure \ref{FIG: 3 Lighthouses - Center - Zoom} we can see the target lighthouse and a way to calculate the dark area behind it. Note that only half of the area behind it is shaded dark gray. Notice the symmetry, $|EB|$, $|FB|$ and $|PB|$ all intersect at the same point behind the target lighthouse. $|EB|$ is the ray coming from top lighthouse, $|FB|$ is the ray coming from bottom lighthouse. This enables us to split the dark area behind the target lighthouse in half. The dark gray shaded area in this case is $d(3)/2$ which is equal to area of the triangle $\triangle ACB$ minus the area of the circular sector formed by $A, C, G$.
\begin{equation*}
Area(ACB) = \frac{x}{2}    
\end{equation*}
The circular sector has an unknown angle $\phi$ but we can overcome this by seeing that $\tan(\phi) = x/1 = x$. This tells us $\phi = \arctan(x)$. Then,
\begin{equation*}
Area(ACG) = \pi \frac{\phi}{360^\circ} = \pi \frac{\arctan(x)}{2\pi} = \frac{\arctan(x)}{2} 
\end{equation*}
\begin{equation}
\label{EQ: d(n)/2 arctan falan}
\frac{d(3)}{2} = Area(ACB) - Area(ACG) = \frac{x - \arctan(x)}{2}
\end{equation}
Remembering that $D(n) = d(n) \times n$ we get $D(3) = d(3) \times 3 = \frac{d(3)}{2} \times 2\times 3$.
\begin{equation}
D(3) = \frac{x - \arctan(x)}{2} \times 2\times 3 = 3(x - \arctan(x))
\end{equation}
Now all that is left to do is find $x$. In this case $120^\circ$ and its complementary $60^\circ$ are beautiful angles, therefore could be useful for us. We already know $|EP| = |PC| = 3$ and $|AC| = 1$ because that is how the problem is defined. First we can find $|EC| = 3\sqrt{3}$ and by Pythagoras rule at the $\triangle EAC$ triangle find $|EA| = \sqrt{26}$. Then, we draw the triangle $\triangle EDP$ and find $|ED|$ to be $3\sqrt{3}/2$. We now have two similar triangles, notice that $\angle CAB = \angle EDB = 90^\circ$, $\angle ABC = \angle EBD$ and $\angle ACB = \angle DEB$. Thanks to this similarity between $\triangle ABC$ and $\triangle EBD$ we can say
\begin{equation*}
\frac{|AB|}{|DB|} = \frac{|BC|}{|EB|} = \frac{|AC|}{|ED|}
\end{equation*}
Looking at $|AC|/|ED| = |BC|/|EB|$ we have
\begin{equation*}
\frac{1}{\frac{3\sqrt{3}}{2}} = \frac{\sqrt{1+x^2}}{\sqrt{26}+x}
\end{equation*}
Squaring both sides
\begin{equation*}
\frac{4}{27} = \frac{1+x^2}{x^2+2x\sqrt{26} + 26}
\end{equation*}
\begin{equation*}
4x^2+8x\sqrt{26}+104 = 27 + 27x^2
\end{equation*}
This gives the equation $23x^2 - 8x\sqrt{26} - 77 = 0$. Using the quadratic formula:
\begin{equation*}
    \begin{split}
         x_{1,2} & = \frac{8\sqrt{26} \pm \sqrt{64\times26 + 4\times23\times77}}{46}
    \end{split}
\end{equation*}
We cannot have negative root because $x$ is the length of $|AB|$, since $8\sqrt{26} < \sqrt{64\times 26 + 4\times 23\times 77}$ we will choose $+$ in place of $\pm$.  Then we have
\begin{equation}
\label{EQ: n=3 X by hand}
x = \frac{\sqrt{1664} + \sqrt{8748}}{46}
\end{equation}
Measuring $x$ using GeoGebra, where we draw these figures, yields $x \approx 2.920$. Calculating \ref{EQ: n=3 X by hand} using a calculator gives $2.92006$. This shows that our calculations are correct.
\begin{equation}
D(3) = 3\left(\frac{\sqrt{1664} + \sqrt{8748}}{46} - \arctan\left(\frac{\sqrt{1664} + \sqrt{8748}}{46}\right)\right)
\end{equation}
The result is $D(3) \approx 5.0376$. 
\subsection{4 \& 6 Lighthouses}
\begin{figure}
    \centering
    \includegraphics[width=0.70\textwidth, frame]{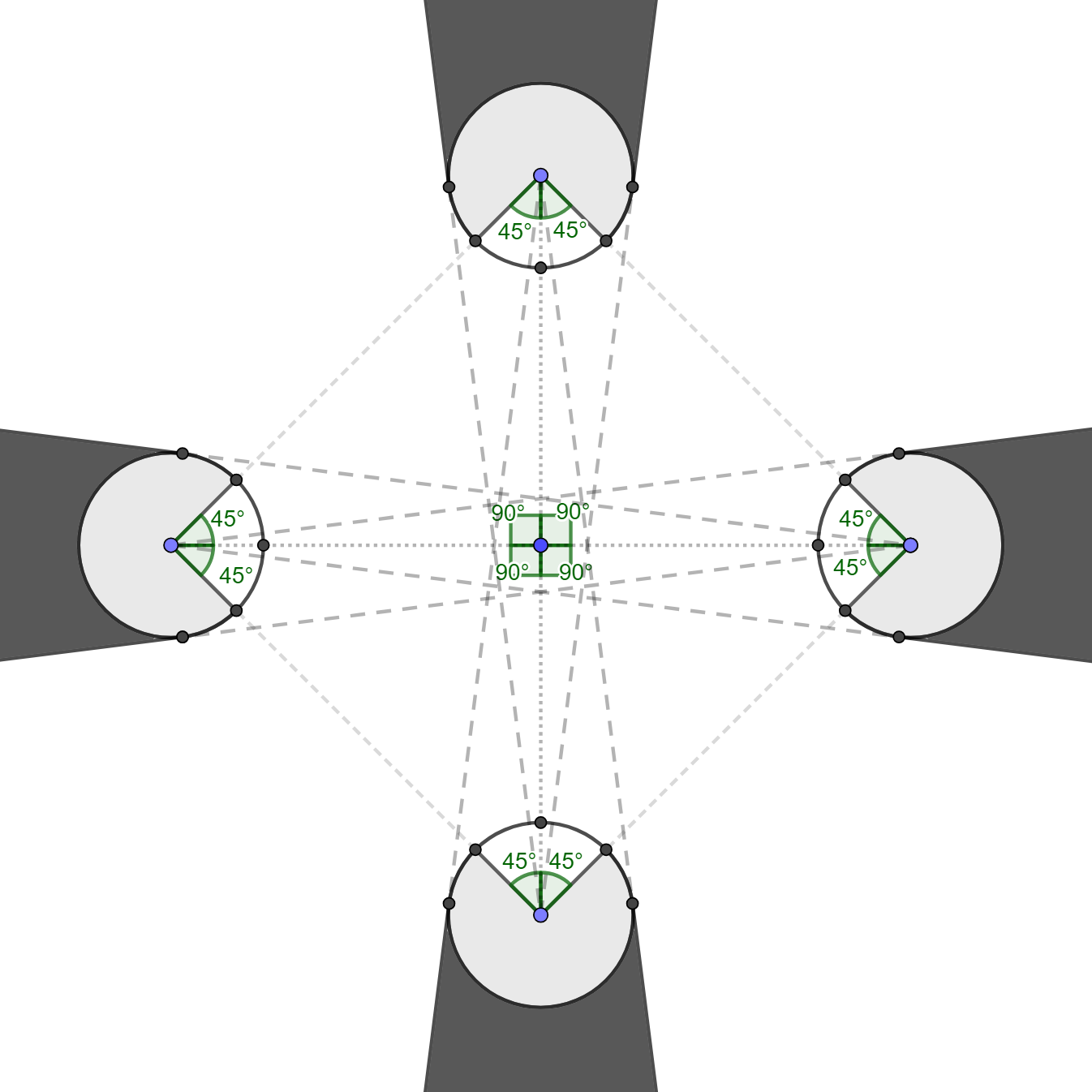} % second figure itself
    \caption{4 Lighthouses.}
    \label{FIG: 4 Lighthouses - Center}
\end{figure}

On figure \ref{FIG: 4 Lighthouses - Center} we can see the case for $n=4$. 
Similar to $n=2$ case we have $D(4) = \infty$ because the dark area behind each lighthouse goes to infinity. It is important to note that this is because no two rays can meet behind a lighthouse. 
\begin{figure}
    \centering
    \includegraphics[width=0.70\textwidth, frame]{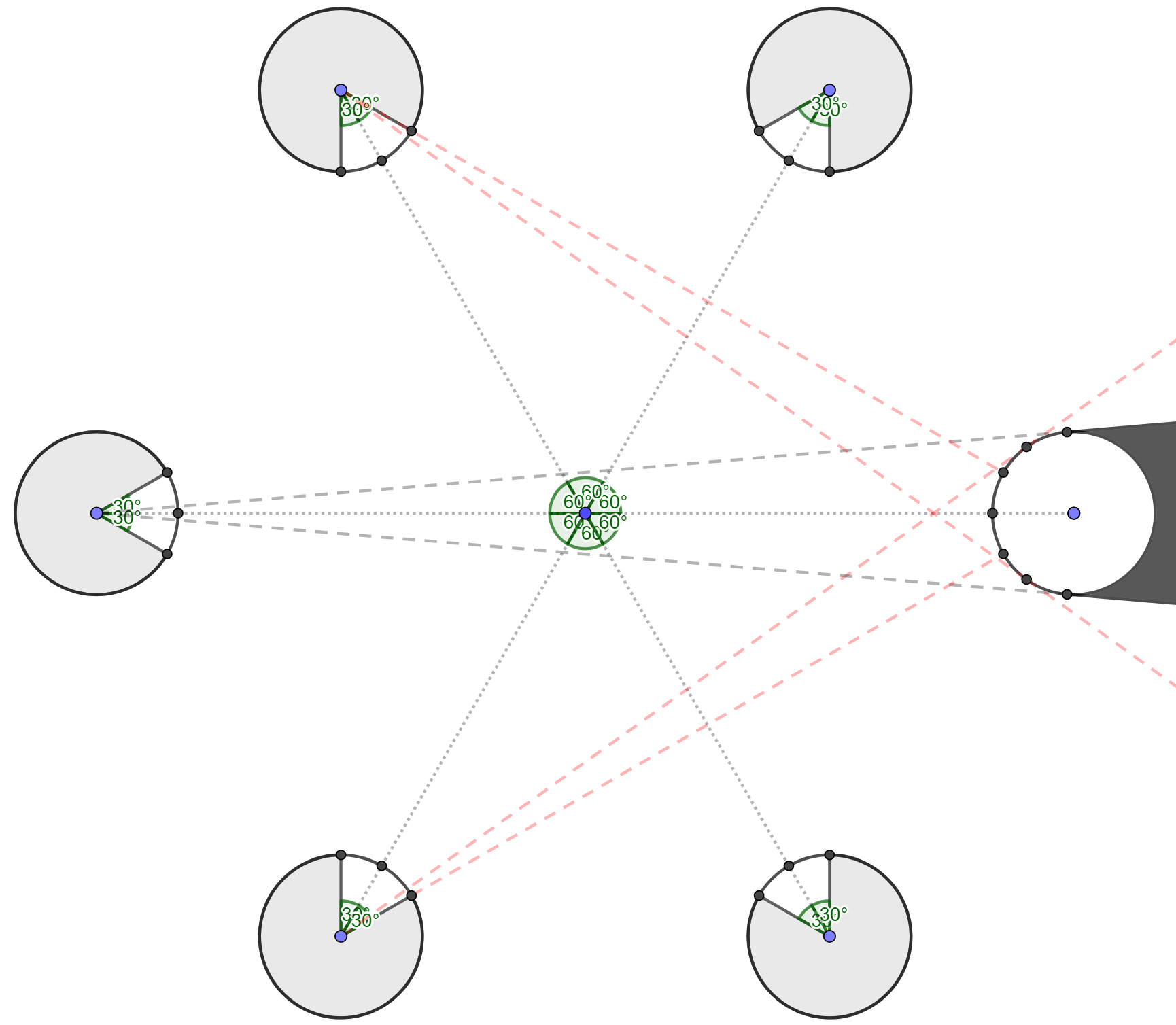}
    \caption{6 Lighthouses and the dark area behind a target lighthouse.}
    \label{FIG: 6 Lighthouses - Center}
\end{figure}
Looking at a target lighthouse for $n=6$ on figure \ref{FIG: 6 Lighthouses - Center} we see the same thing happening, $D(6) = \infty$
\subsection{Even Number of Lighthouses}
For even lighthouses, we were not able to draw a tangent to define a finite dark area behind the target lighthouse. This causes the dark area to be infinite. We can show that this is always true, using proof by contradiction.
\begin{theorem}
\label{THEO: Center - Even - Theorem}
$D(n) = \infty, n \equiv 0 \pmod{2}$
\end{theorem}
\begin{proof}
Imagine an even number of lighthouses placed around a center. Place the first lighthouse directly $n$ units to the right of the placement center (as we always did in our figures). Let $L_0$ be the first and target lighthouse, number the rest of the lighthouses as $L_1, L_2, ..., L_{n-1}$ counter-clockwise. If we number like this, the lighthouse on the opposite side of $L_0$ is $L_{n/2}$. We have to show that none of the lighthouses numbered $L_{n/2 - 1}, L_{n/2 - 2}, ..., L_1$ are able to define a finite dark area behind $L_0$. Note how they are all numbered in reference to $L_{n/2}$. Now to use proof by contradiction, we assume that it is possible to draw a finite dark area defining tangent from one of the $L_{n/2 - 1}, L_{n/2 - 2}, ..., L_1$ giving us the figure \ref{FIG: Center Even Proof}. In this figure, $E$ is the center of lighthouse $L_{n/2-k}, k \in \mathbb{N}$, $P$ is the placement center and $C$ is the center of target lighthouse $L_0$. By definition $|EP| = |PC|$, but if $k=1$ then $\angle EBP = \angle PEB$ which requires $|EP| = |PB|$. $|PB| \ne |PC|$ therefore $k=1$ is not possible. For $k>1$ we realize $\angle EBP > \angle PEB$. Now looking at $|EP| = n$ and $|PB| = n + |CB|$, if $\angle EBP > \angle PEB$ this would require $|EP| > |PB|$,
\begin{equation}
\label{EQ: Center Even Proof eq}
\begin{split}
    |EP| &> |PB| \\
    n &> n + |CB| \\
    0 &> |CB|
\end{split}
\end{equation}
Thus showing that this is not possible. Another case to consider is when there is such a tangent that $\angle PEB < \alpha/2$, but again this causes $\angle EBP$ to be bigger and inevitably  $\angle EBP > \angle PEB$, which we have shown to be contradictory in \ref{EQ: Center Even Proof eq}.
\begin{figure}
    \centering
    \includegraphics[width=0.95\textwidth]{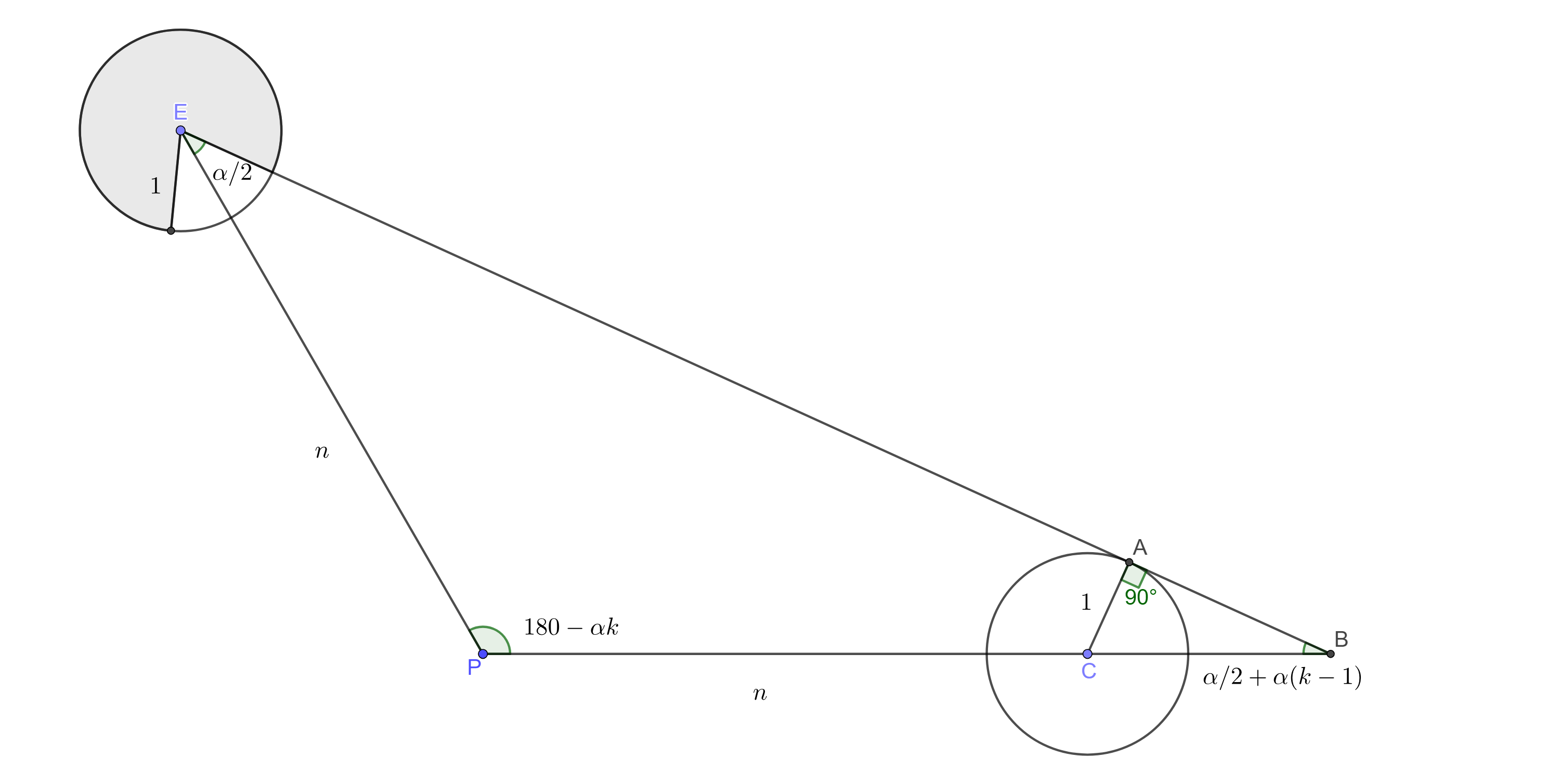}
    \caption{$D(n) = \infty, n \equiv 0 \pmod{2}$ proof by contradiction.}
    \label{FIG: Center Even Proof}
\end{figure}
\end{proof}
Now we can safely say that dark area is infinite for even number of lighthouses.

\subsection{Odd Number of Lighthouses}
\begin{figure}
    \centering
    \includegraphics[width=0.75\textwidth, frame]{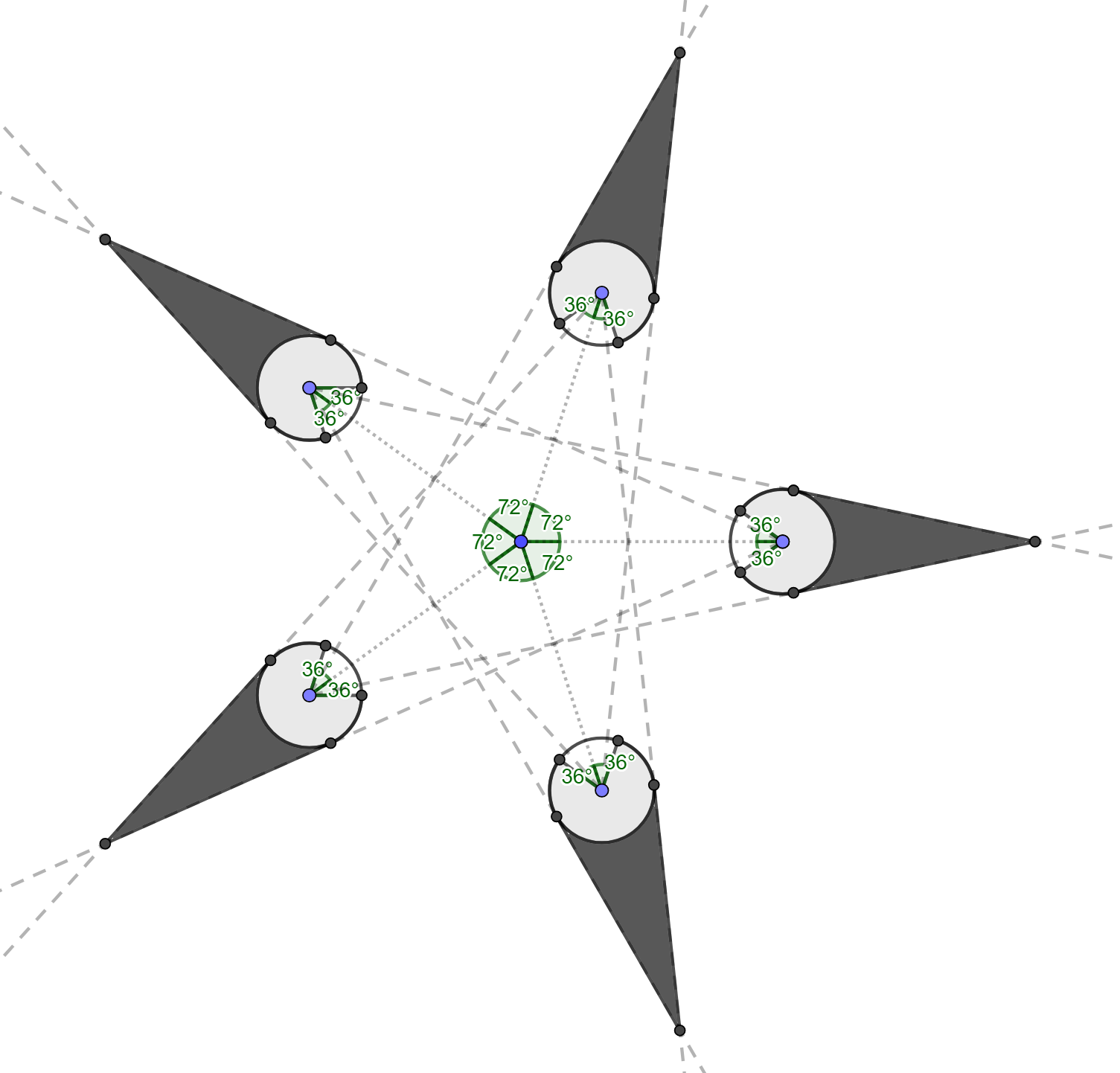}
    \caption{5 Lighthouses.}
    \label{FIG: 5 Lighthouses - Center}
\end{figure}

Looking at 5 lighthouses case in figure \ref{FIG: 5 Lighthouses - Center} it seems as if the furthest lighthouses define the dark area behind the target lighthouse.
We will show that it is indeed the furthest lighthouses that define the dark area behind a target lighthouse and then give a formula to find the dark area itself.
\begin{theorem}
The finite dark area $d(n)$ behind a target lighthouse is defined by the tangents coming from two furthest lighthouses when $n \equiv 1 \pmod{2}$.
\end{theorem}
\begin{proof}
\label{PROOF: Center - Even - Theorem}
Our proof is similar to the proof we had for theorem \ref{THEO:  Center - Even - Theorem}. Imagine an odd number of $n$ lighthouses. Place the first lighthouse directly $n$ units to the right of the placement center. Let $L_0$ be the first and target lighthouse. Draw a line that passes through the placement center and the center of the target lighthouse. This divides the plane in half, $(n-1)/2$ lighthouses on one side and $(n-1)/2$ on the other. Since $n$ is an odd number we can say $n = 2m + 1$. So in other words, we have $m$ lighthouses on one side and $m$ on the other. Starting from the right, number the lighthouses on the upper half as $L_1, L_2, ..., L_m$, counter-clockwise. Our claim is that $L_0$ is illuminated by $L_m$ and $L_m$ only. We said ``furthest lighthouses'' which is plural, the other lighthouses is the lighthouse that is $L_m$'s reflection on the line we just drew. Like we demonstrated on $n=3$ case the dark area behind the lighthouse is divided in two by this line, calculating on one half suffice. We can use proof by contradiction, assuming that some other lighthouse $L_i, 0<i<m$ illuminates $L_0$. Drawing the figure \ref{FIG: Center Odd Proof} gives us an idea.
\begin{figure}[h!]
    \centering
    \includegraphics[width=0.95\textwidth]{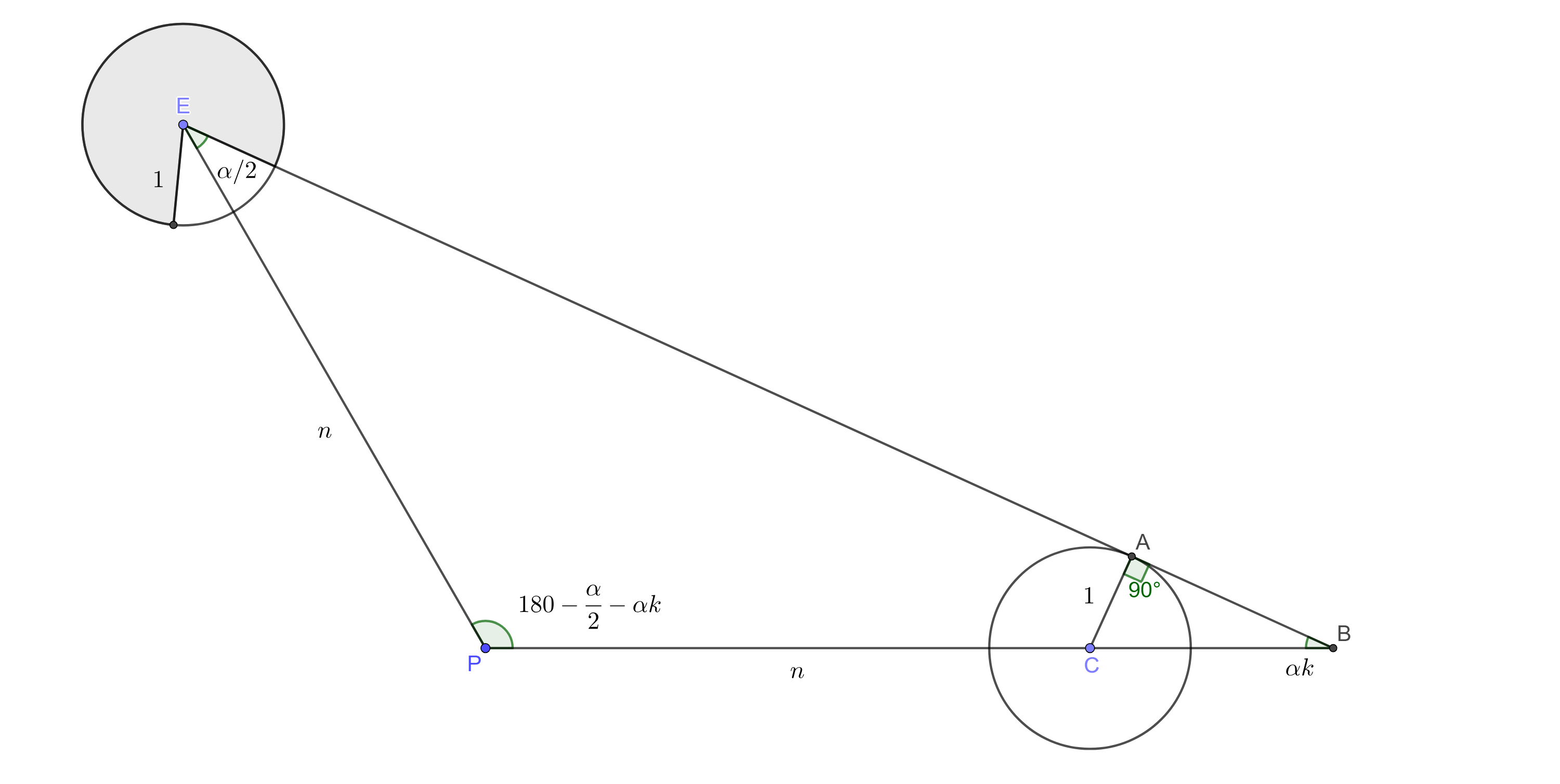}
    \caption{$L_0$ and the two furthest lighthouses.}
    \label{FIG: Center Odd Proof}
\end{figure}
In this figure, $E$ is the center of lighthouse $L_{m-k}, k \in \mathbb{N} \cup \{0\}$, $P$ is the placement center, $C$ is the center of the target lighthouse $L_0$. For positive values of $k$ we see that $\alpha k > \alpha/2$. This would require $|EP| > |PB|$ but that is not possible, exactly the same contradiction on \eqref{EQ: Center Even Proof eq}. We now know that $k=0$, so $\angle EPB = 180^\circ - \alpha/2$ and we are illuminated by $L_m$, but what is $\angle EPB$? Let us say that $\angle EBP = \beta$ and $PEB = \alpha/2 - \beta$ such that $\angle EPB < \angle PEB$. This gives us $\beta < \alpha/2 - \beta$ which results in $\beta < \alpha/4$. As the number of lighthouses increase, the angle $\angle EPB$ approaches $180^\circ$ because $\alpha$ gets smaller, and now we also see that $\beta$ gets smaller too, such that $\beta < \alpha/4$. In conclusion, $L_m$ does illuminate $L_0$ and it is the only one doing so.
\end{proof}
\begin{theorem}
$D(n) = n(x_n - \arctan(x_n)), n \equiv 1 \pmod{2}, n > 1$ where
\begin{equation}
\label{EQ: x_n in theorem for center}
    x_n = \frac{\sqrt{4n^2\cos^2\left(\frac{\pi}{2n}\right) - 1} + 2n^2\sin(\frac{\pi}{n})\cos^2(\frac{\pi}{2n})}{ n^2\sin^2\left(\frac{\pi}{n}\right)-1}
\end{equation}
\end{theorem}
\begin{proof}
This is going to be a pretty straightforward proof. To start, let us again imagine the same figure we did for the previous proof \ref{PROOF: Center - Even - Theorem}. We have $L_0$ to the right and $L_1, L_2, ..., L_m$ going counter-clockwise. We just showed that we only need to care about $L_m$ so that is exactly what we are going to do.
\begin{figure}[h!]
    \centering
    \includegraphics[width=1.2\textwidth]{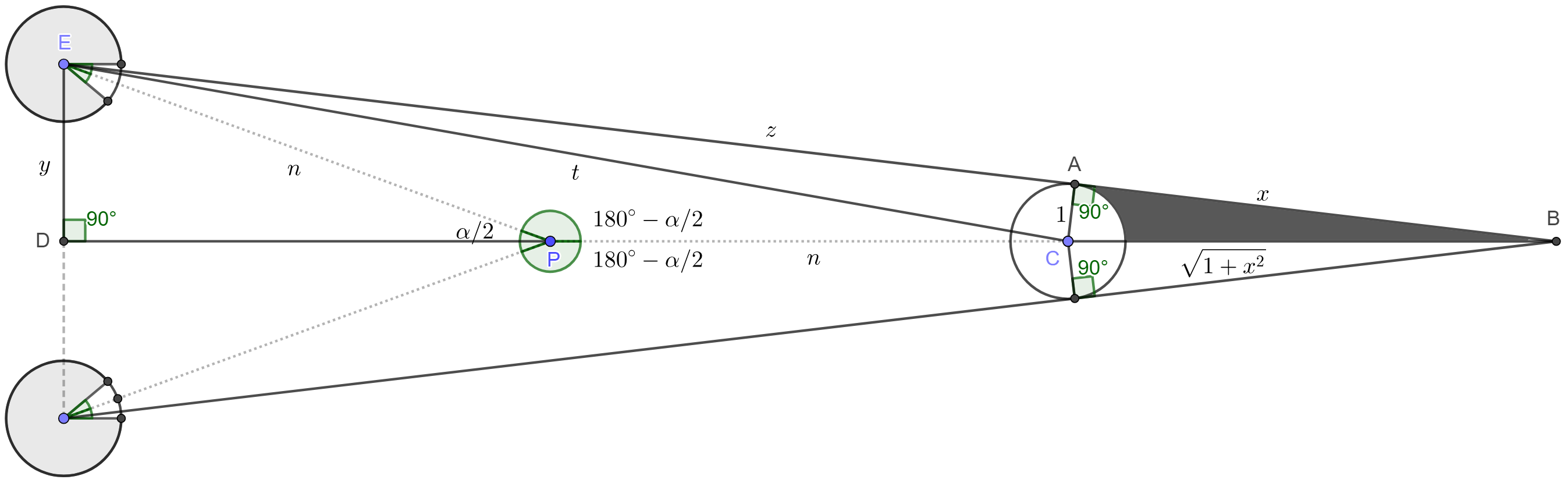} 
    \caption{$n$ Lighthouses, $n \equiv 1 \pmod{2}$, focusing on $L_0$, the lighthouse to the upper-left is $L_m$}
    \label{FIG: N Odd Lighthouses - Center}
\end{figure}
We get figure \ref{FIG: N Odd Lighthouses - Center} as a result. We have some new variables in our figure, namely $y, z$ and $t$. We will describe all three of them in terms of $n$ and then work our way towards $x$. $x$ in this figure is $x_n$ in equation \eqref{EQ: x_n in theorem for center}. Like we have done for $n=3$, if we find $x$ we can find $d(n)$, and consequently $D(n)$. So, let us describe the uninvited guests one by one:
\begin{itemize}
    \item $y$ can be found by looking at $\sin(\alpha/2) = y/n$, which gives
    \begin{equation}
    y = n \times \sin\left(\frac{360^\circ}{2n}\right) = n \times \sin\left(\frac{180^\circ}{n}\right)
    \end{equation}
    \item $t$ can be found using Cosine Theorem on $|EC|$ seen by the angle $\angle EPC$. The theorem gives us $t^2 = n^2 + n^2 - 2n^2\cos(180^\circ - \alpha/2)$. The expression $n^2 + n^2 - 2n^2\cos(180^\circ - \alpha/2)$ can be reduced.
    \begin{equation*}
     n^2 + n^2 - 2n^2\cos(180^\circ - \alpha/2) = 2n^2(1-\cos(180^\circ - \alpha/2)) 
    \end{equation*}
    Since $1-\cos(\phi) = 2\sin ^2 (\phi/2)$,
    \begin{equation*}
        \begin{split}
           2n^2(1-\cos(180^\circ - \alpha/2)) & = 2n^2(2\sin^2((180^\circ - \alpha/2)/2)) \\
                & = 4n^2\left(\sin^2\left(\frac{180^\circ - \frac{\alpha}{2}}{2}\right)\right) \\
                & = 4n^2\left(\sin^2\left(\frac{180^\circ - \frac{360^\circ}{2n}}{2}\right)\right) \\
                & = 4n^2\left(\sin^2\left(90^\circ \left(1 - \frac{1}{n}\right)\right)\right)
        \end{split}
    \end{equation*}
    Finally we get
    \begin{equation}
    t = \sqrt{4n^2\left(\sin^2\left(90^\circ \left(1 - \frac{1}{n}\right)\right)\right)} = 2n \times \sin\left(90^\circ \left(1 - \frac{1}{n}\right)\right)    
    \end{equation}
    \item $z$ can be written in terms of $t$ using Pythagoras rule at $\triangle EAC$. The rule gives us $t^2 = 1 + z^2$ therefore $z^2 = t^2 - 1$. Since we already wrote $t$ in terms of $n$ we get
    \begin{equation}
    z^2 =  4n^2\left(\sin^2\left(90^\circ \left(1 - \frac{1}{n}\right)\right)\right) - 1    
    \end{equation}
\end{itemize}
Now that we have $y, z$ and $t$ in terms of $n$ we can work on $x$. We will be doing the same thing we did for $n=3$, thanks to the similarity $\triangle EDB \sim \triangle CAB$.
\begin{equation*}
\frac{|AB|}{|DB|} = \frac{|BC|}{|EB|} = \frac{|AC|}{|ED|}    
\end{equation*}
Looking at $|AC|/|ED| = |BC|/|EB|$ we have
\begin{equation*}
    \frac{1}{y} = \frac{\sqrt{1+x^2}}{z + x}
\end{equation*}
Squaring both sides
\begin{equation*}
    \frac{1}{y^2} = \frac{1+x^2}{x^2 + 2xz + z^2}
\end{equation*}
This gives us $x^2 + 2xz + z^2 = y^2 + x^2y^2$ which yields the equation
\begin{equation*}
    x^2(y^2-1) + x(-2z) + (y^2 - z^2) = 0
\end{equation*}
Plugging this into quadratic formula gives us
\begin{equation*}
    x_{1,2} = \frac{2z \pm \sqrt{4z^2  - 4(y^2-1)(y^2-z^2)}}{2y^2 - 2}
\end{equation*}
By definition of the problem, $y < z$ and because of this the $\pm$ will have to be $+$ sign, otherwise $x$ will be negative and $|AB|$ can't be negative.
\begin{equation*}
    \begin{split}
        x &= \frac{2z + \sqrt{4z^2  - 4(y^2-1)(y^2-z^2)}}{2y^2 - 2} \\
        & = \frac{2z + \sqrt{4z^2  - 4y^4 + 4y^2z^2 + 4y^2 - 4z^2}}{2y^2 - 2} \\
        & = \frac{2z + \sqrt{4y^2(z^2-y^2+1)}}{2y^2 - 2} \\
        & = \frac{2z + 2y\sqrt{z^2-y^2+1}}{2y^2 - 2} \\
        & = \frac{z + y\sqrt{z^2-y^2+1}}{y^2-1}
    \end{split}
\end{equation*}
We should take a moment and write what $y^2$ and $z^2$ are using radians instead of degrees.
\begin{equation*}
    y^2 =  n^2\sin^2\left(\frac{\pi}{n}\right)
\end{equation*}
\begin{equation*}
    z^2 = 4n^2\sin^2\left(\frac{\pi}{2} - \frac{\pi}{2n}\right) - 1 = 4n^2\cos^2\left(\frac{\pi}{2n}\right) - 1
\end{equation*}
We can see that $\sqrt{z^2-y^2+1}$ is
\begin{equation*}
    \begin{split}
       \sqrt{z^2-y^2+1} &= \sqrt{4n^2\sin^2\left(\frac{\pi}{2} - \frac{\pi}{2n}\right) - 1 - n^2\sin^2\left(\frac{\pi}{n}\right) + 1} \\
       &= \sqrt{n^2\left(4\sin^2(\pi/2-\pi/2n) - \sin^2(\pi/n)\right)} \\
       &= \sqrt{n^2\left(4\cos^2(\pi/2n) - \sin^2(\pi/n)\right)}
    \end{split}
\end{equation*}
Remember that $\sin{2\phi} = 2\sin{\phi}\cos{\phi}$. If $\phi = \pi/2n$ we get $\sin(\pi/n) = 2\sin(\pi/2n)\cos(\pi/2n)$ therefore $\sin^2(\pi/n) = 4\sin^2(\pi/2n)\cos^2(\pi/2n)$. We further reduce:
\begin{equation*}
    \begin{split}
       \sqrt{z^2-y^2+1} & = \sqrt{n^2(4\cos^2(\pi/2n)-4\sin^2(\pi/2n)\cos^2(\pi/2n))}\\
       &= \sqrt{4n^2\cos^2(\pi/2n)(1-\sin^2(\pi/2n))}
    \end{split}
\end{equation*}
Also remembering $\sin^2\phi + cos^2\phi = 1$ we can further reduce
\begin{equation*}
    \begin{split}
       \sqrt{z^2-y^2+1} & = \sqrt{4n^2\cos^2(\pi/2n)(\cos^2(\pi/2n))} \\
       & = 2n\cos^2(\pi/2n)
    \end{split}
\end{equation*}
$y\sqrt{z^2-y^2+1}$ is then $2n^2\sin(\pi/n)\cos^2(\pi/2n)$. Finally we have the formula for $x$:
\begin{equation}
\label{EQ: Center X by formula}
    x = \frac{\sqrt{4n^2\cos^2\left(\frac{\pi}{2n}\right) - 1} + 2n^2\sin(\frac{\pi}{n})\cos^2(\frac{\pi}{2n})}{ n^2\sin^2\left(\frac{\pi}{n}\right)-1}
\end{equation}
\end{proof}
Giving $n=3$ in equation \eqref{EQ: Center X by formula} and using calculator yields 
\begin{equation*}
x = \frac{\sqrt{416} + \sqrt{2187}}{23} = \frac{\sqrt{1664}+\sqrt{8748}}{46}    
\end{equation*}
which is exactly same as the result we obtained by hand at equation \eqref{EQ: n=3 X by hand}. Looking back at $n=5$ case, measuring $x$ in GeoGebra yields $x \approx 4.7190$. Plugging $n=5$ to equation \eqref{EQ: Center X by formula} returns $x \approx 4.7190$. When it comes to $d(n)$ we can do what we did in \eqref{EQ: d(n)/2 arctan falan} and by looking at figure \ref{FIG: N Odd Lighthouses - Center} we can say that $d(n) = x - \arctan{x}$. Taking $x=4.7190$ and plugging it in $d(5)$ yields
\begin{equation*}
    d(5) = 4.7190 - \arctan{4.7190} \approx 3.3570
\end{equation*}
Multiplying this by 5 and we get
\begin{equation}
    D(5) = d(5) \times 5 \approx 16.7851
\end{equation}
%% END %%

%% file: Sections/2.3-Arc-Light-Source.tex
%% BEGIN %%
\section{Point Light Sources at the Arc}
The second variation of the problem is when the points on the arc seen by the illumination angle $\alpha$ act as point light sources. This case is a bit more complex than the ``Point Light Source at the Center'' variation, which we have solved. We will be giving examples for 1, 2, 3, 4 and 5 lighthouses, then try to generalize it, which is where the problem passionately slaps us for trying to do so.
\subsection{1 Lighthouse}
This is the same case as it was for point light source at the center, because the illumination angle is $360^\circ$. To see how this looks like, refer to figure \ref{FIG: 1 Lighthouse}. For the record, $D(1) = 0$.
\subsection{2 Lighthouses}
\begin{figure}
    \centering
    \includegraphics[width=0.70\textwidth, frame]{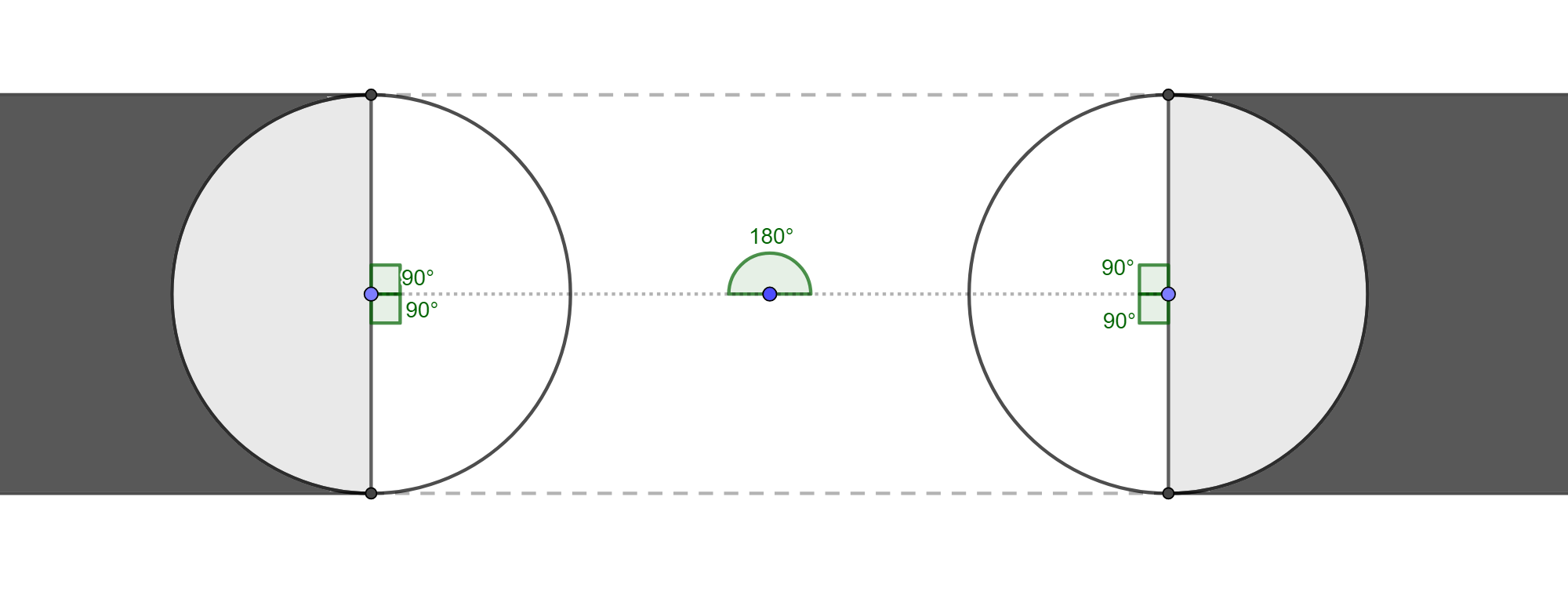}
    \caption{2 Lighthouses}
    \label{FIG: 2 Lighthouses - Arc}
\end{figure}

This case is shown on figure \ref{FIG: 2 Lighthouses - Arc}.  
Again, the dark area is infinite, but in a slightly different way. It is literally ``barely'' infinite. The point light sources $180^\circ$ apart from each other in both arcs make a parallel light ray tangent to the lighthouses facing each other, resulting in a rectangular-like dark area that extends to infinity. We are safe to say $D(2) = \infty$.
\subsection{3 Lighthouses}
\begin{figure}
    \centering
    \includegraphics[width=0.70\textwidth, frame]{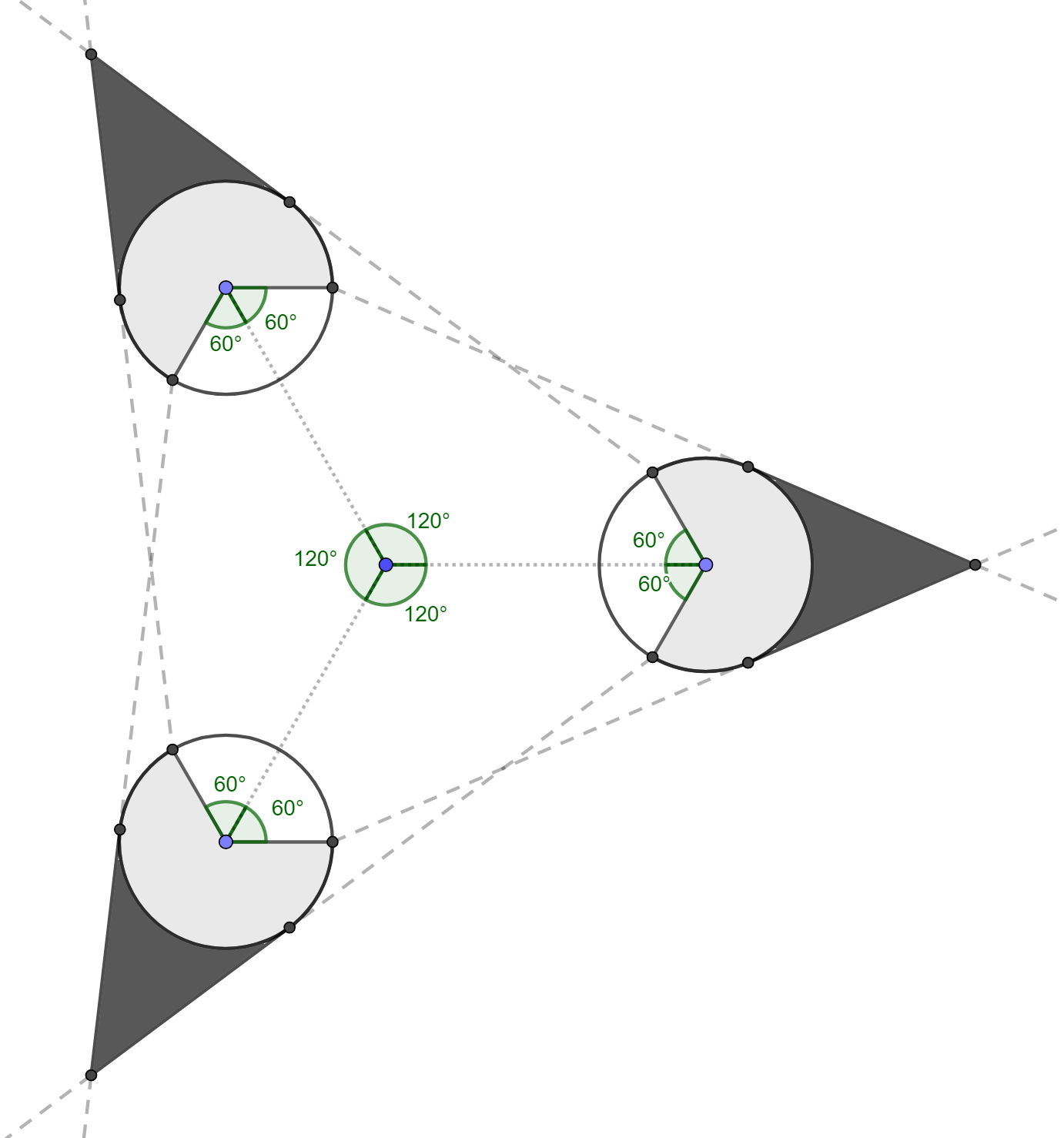}
    \caption{3 Lighthouses}
    \label{FIG: 3 Lighthouses - Arc}
\end{figure}

\begin{figure}
    \centering
    \includegraphics[width=0.70\textwidth, frame]{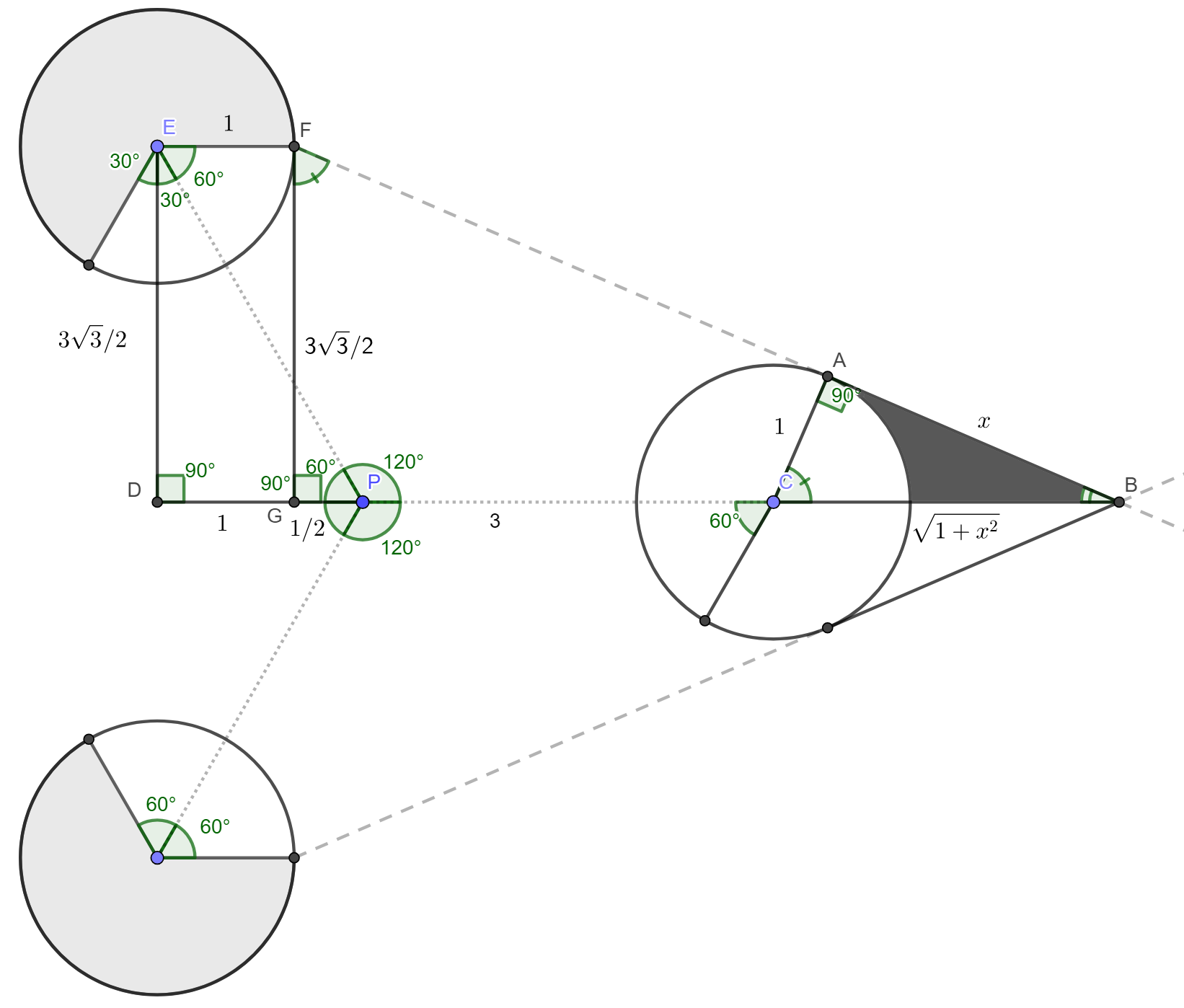}
    \caption{3 Lighthouses with the target lighthouse}
    \label{FIG: 3 Lighthouses - Arc - Zoom}
\end{figure}

Figure \ref{FIG: 3 Lighthouses - Arc} shows us the case for $n=3$. 
Like we always did so far we will be focusing on a single lighthouse. On figure \ref{FIG: 3 Lighthouses - Arc - Zoom} we can see a way to approach $x$.
We should immediately warn that though it looks as if $\angle ACB = 60^\circ$ that is not correct, it is just that the angle is quite close to $60^\circ$. To find $x$ we draw the points $D$ and $G$, thereby finding $|DG| = 1$ and $|GP| = 1/2$. We already know $|PC| = 3$. Since the $30^\circ$ angle $\angle DEP$ is seeing $|DP| = 3/2$ then $60^\circ$ angle $\angle EPD$ sees $|ED| = 3\sqrt{3}/2$. Also because $|ED| = |FG|$ we have $|FG| = 3\sqrt{3}/2$. We the notice the similarity $\triangle FGB \sim \triangle CAB$ which tells us that there is a ratio $|AC|/|FG| = |AB|/|GB|$. This gives us the equation below.
\begin{equation*}
    \frac{1}{\frac{3\sqrt{3}}{2}} = \frac{x}{\frac{7}{2} + \sqrt{1+x^2}}
\end{equation*}
The author could not reduce this to a quadratic equation and sought help from Wolfram$|$Alpha \cite{wolfram-n3-edge}, which gave $x =  3(4\sqrt{2} + 7\sqrt{3})/23 \approx 2.3192$. We then used GeoGebra to measure the length $|AB|$ which gives approximately $2.3192$ so we can assume that the equation is correct. Taking $x=2.3192$ the total dark area is 
\begin{equation}
    D(3) = 3(2.3192 - \arctan(2.3192)) \approx 3.4665
\end{equation}
\subsection{4 Lighthouses}
\begin{figure}
    \centering
    \includegraphics[width=0.9\textwidth, frame]{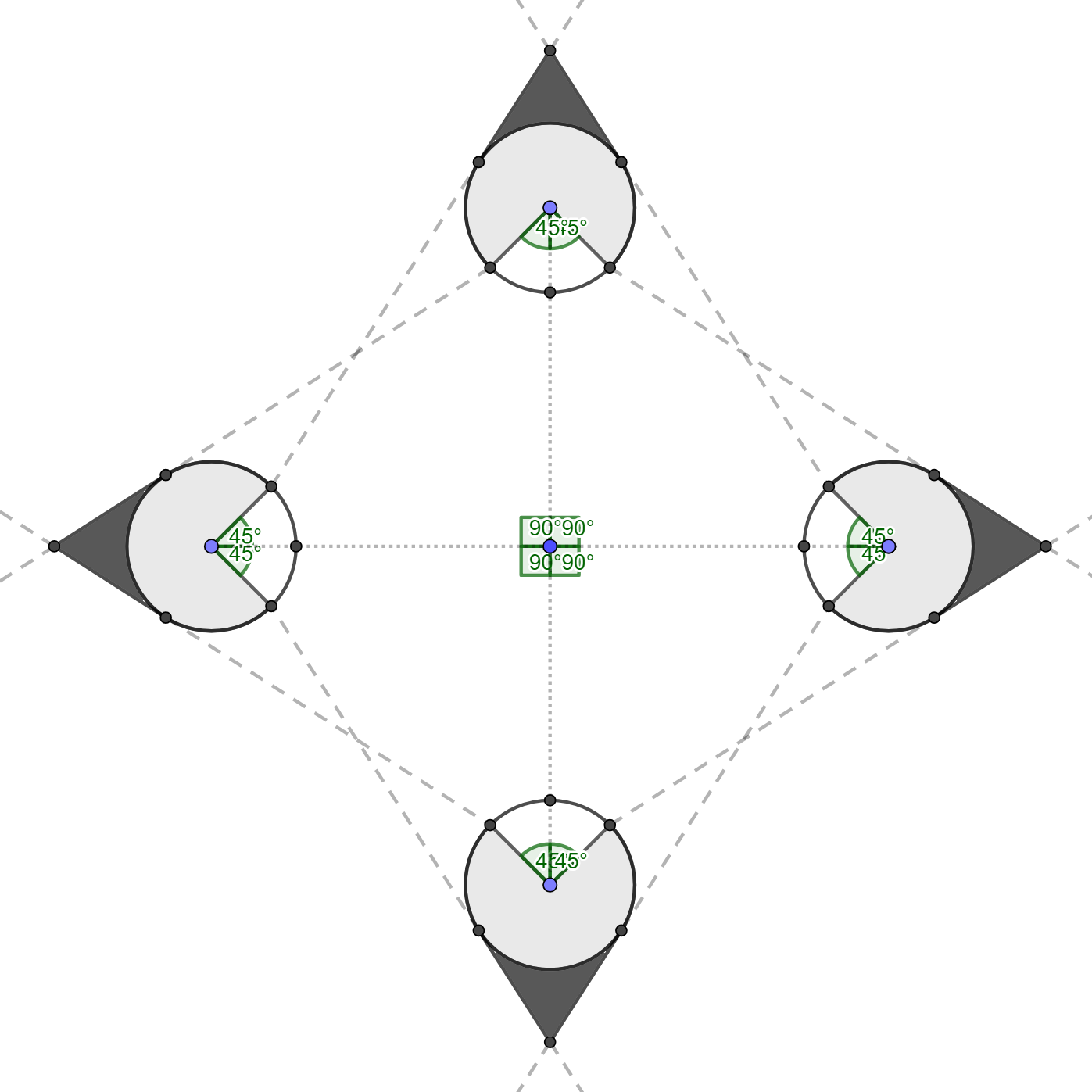} % second figure itself
    \caption{4 Lighthouses, Arc of Point Light Sources.}
    \label{FIG: 4 Lighthouses - Edge}
\end{figure}

\begin{figure}
    \centering
    \includegraphics[width=\textwidth]{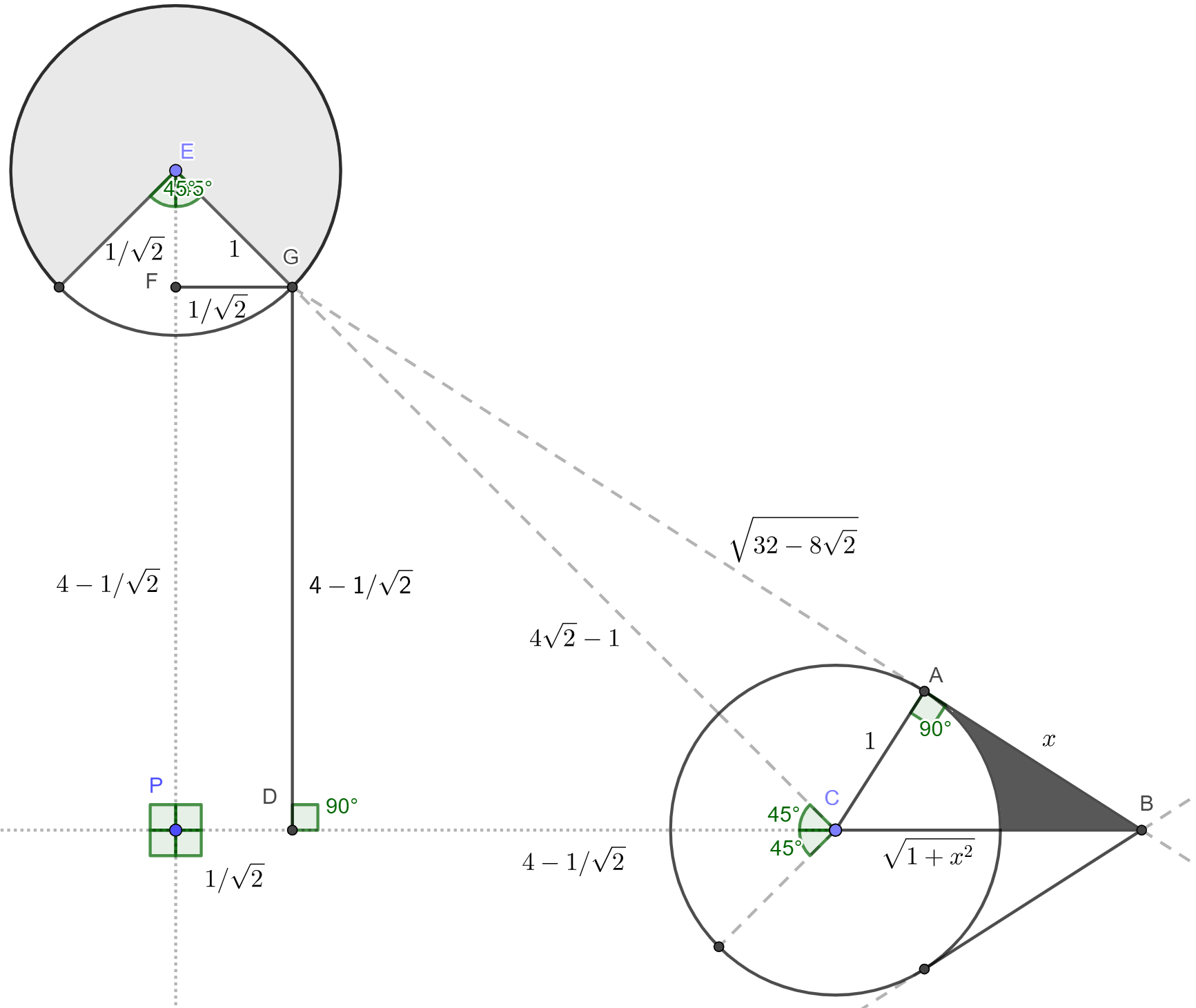} % second figure itself
    \caption{4 Lighthouses, focusing on a single lighthouse.}
    \label{FIG: 4 Lighthouses - Edge - Focusing}
\end{figure}

Unlike the first variation, we actually have a finite dark area in this case. Figure \ref{FIG: 4 Lighthouses - Edge} shows the case for $n=4$. 
Now, the closest neighbors of a target lighthouse can illuminate behind it. In figure \ref{FIG: 4 Lighthouses - Edge - Focusing} we will try to find $x$.
We find $|FG| = 1/\sqrt{2}$ with a simple glance at $\triangle EFG$. Thanks to the right-angled triangle $\triangle GDC$ where $|GD| = |DC|$ we can use Pythagoras rule to find $|GC| = 4\sqrt{2}-1$. Once again using Pythagoras rule we find $|GA| = \sqrt{32 - 8\sqrt{2}}$. A final Pythagoras rule $|GD|^2 + |DB|^2 = |GB|^2$ gives us
\begin{equation*}
    \left(4 - \frac{1}{\sqrt{2}}\right)^2 + \left(4- \frac{1}{\sqrt{2}} + \sqrt{1+x^2}\right)^2 = \left(\sqrt{32 - 8\sqrt{2}} + x\right)^2
\end{equation*}
After unwrapping and reducing both sides we get
\begin{equation*}
    1 + \left(4 -  \frac{1}{\sqrt{2}}\right)\sqrt{1 + x^2} = x\sqrt{32-8\sqrt{2}}
\end{equation*}
Again, the author was unable to reduce and solve this. We use Wolfram$|$Alpha  \cite{wolfram-n4-edge} to find $x$ which returns $x \approx 1.5637$. Using GeoGebra to measure the length $|AB|$ also yields approximately $1.5637$. Taking $x=1.5637$ we can find $D(4)$.
\begin{equation}
    D(4) = 4(1.5637 - \arctan(1.5637)) \approx 2.24745
\end{equation}
\subsection{5 Lighthouses}
\begin{figure}
    \centering
    \includegraphics[width=0.9\textwidth, frame]{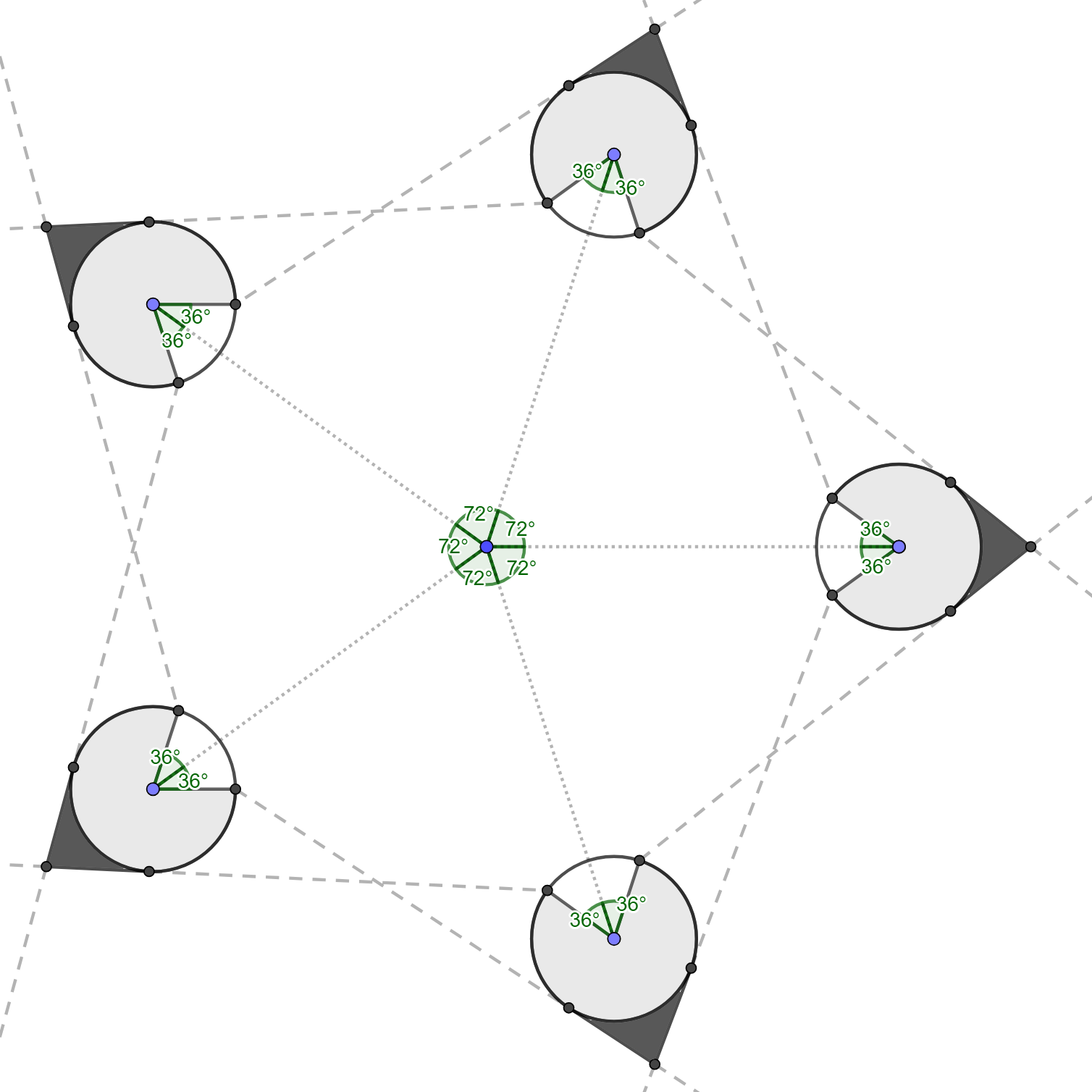} % second figure itself
    \caption{5 Lighthouses}
    \label{FIG: 5 Lighthouses - Edge}
\end{figure}

\begin{figure}[ht]
    \centering
    \includegraphics[width=\textwidth]{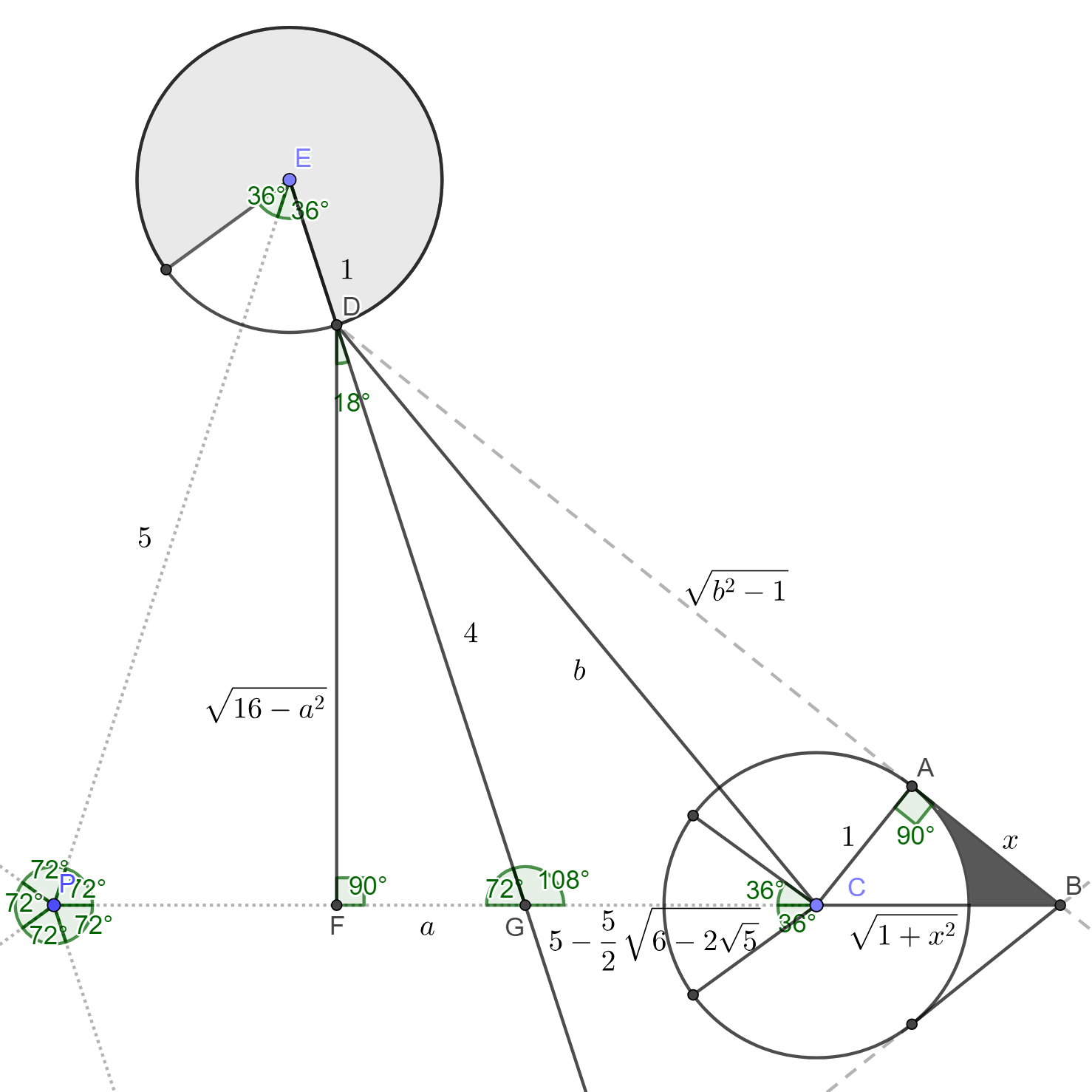} % second figure itself
    \caption{5 Lighthouses, focusing on a single lighthouse.}
    \label{FIG: 5 Lighthouses - Edge - Focused}
\end{figure}

Looking at $n=5$ on figure \ref{FIG: 5 Lighthouses - Edge}, again it appears that the closest neighbors define the dark area behind a target lighthouse. 
Looking closely to the target lighthouse as shown in figure \ref{FIG: 5 Lighthouses - Edge - Focused}
Immediately we are blessed with $|EP| = |EG| = 5$, therefore $|DG| = 4$. We want to use the similarity $\triangle DBF \sim \triangle ABC$ like we did on previous cases, so we have to find $\sqrt{b^2-1}$ and $\sqrt{16-a^2}$ before we actually come to $x$. To find $a$ we can use the Cosine Theorem as $|FG|^2 = |DF|^2 + |DG|^2 - 2|DF||DG|\cos(18^\circ)$ and leave $a$ out to find its value. Remembering that $\cos(18^\circ) = \sqrt{10 + 2\sqrt{5}}/4$ we get
\begin{equation*}
    a^2 = 16 + 16 - a^2 - 2(4)(\sqrt{16 - a^2})\frac{\sqrt{10+2\sqrt{5}}}{4}
\end{equation*}
This reduces to
\begin{equation*}
    2(16 - a^2) = 2(\sqrt{16 - a^2})\sqrt{10+2\sqrt{5}}
\end{equation*}
\begin{equation*}
    \frac{16 - a^2}{\sqrt{16 - a^2}} = \sqrt{10+2\sqrt{5}}
\end{equation*}
\begin{equation*}
    \sqrt{16 - a^2} = \sqrt{10+2\sqrt{5}}
\end{equation*}
Squaring both sides,
\begin{equation*}
    16 - a^2 = 10+2\sqrt{5}
\end{equation*}
\begin{equation*}
    a^2 = 6 - 2\sqrt{5}
\end{equation*}
\begin{equation}
    a = \sqrt{6 - 2\sqrt{5}}
\end{equation}
Next, we find $|GC|$. $|GC| = 5 - |PG|$ and $|PG|$ can be found by Cosine Theorem as $|PG|^2 = 5^2 + 5^2 - 2(5^2)(5^2)\cos(36^\circ)$. This is reduced to $|PG| = 50(1 - \cos(36^\circ)) = 50(1 - (1+\sqrt{5})/4)$.
\begin{equation*}
    |PG|^2 = \frac{50}{4}(3 - \sqrt{5})
\end{equation*}
\begin{equation*}
    |PG| = \frac{5}{2}\sqrt{6 - 2\sqrt{5}}
\end{equation*}
\begin{equation*}
    |GC| = 5 - \frac{5}{2}\sqrt{6 - 2\sqrt{5}}
\end{equation*}
Now that we found $|GC|$ and $a$ we can find $|FC| = a + |GC|$.
\begin{equation}
    |FC| = 5 - \frac{5}{2}\sqrt{6 - 2\sqrt{5}} + \sqrt{6 - 2\sqrt{5}} = 5 - \frac{3}{2}\sqrt{6 + 2\sqrt{5}}
\end{equation}
Next, we can work on $b$ by doing Pythagoras rule as $b^2 = |FC|^2 + |DF|^2$.
\begin{equation*}
    b^2 = (5 - \frac{3}{2}\sqrt{6 - 2\sqrt{5}})^2 + 16 - a^2
\end{equation*}
\begin{equation*}
    b^2 = 25 + \frac{9}{4}(6 - 2\sqrt{5}) -15\sqrt{6 - 2\sqrt{5}} + 16 - 6 + 2\sqrt{5}
\end{equation*}
\begin{equation*}
    b^2 = 41 + \frac{5}{4}(6 - 2\sqrt{5}) -15\sqrt{6 - 2\sqrt{5}} 
\end{equation*}
Looking at $|DA| = \sqrt{b^2-1}$ we find 
\begin{equation*}
     \sqrt{b^2-1} = \sqrt{40 + \frac{5}{4}(6 - 2\sqrt{5}) -15\sqrt{6 - 2\sqrt{5}} }
\end{equation*}
We have everything ready to find $x$. We will do it by using the $\triangle DBF \sim \triangle ABC $. The similarity gives us the equation $|AC|/|DF| = |CB|/|DB|$.
\begin{equation*}
    \frac{1}{\sqrt{16-a^2}} = \frac{\sqrt{1+x^2}}{x + \sqrt{b^2-1}}
\end{equation*}
Squaring both sides
\begin{equation*}
    \frac{1}{16-a^2} = \frac{1+x^2}{x^2 + b^2-1 + 2x \sqrt{b^2-1}}
\end{equation*}
\begin{equation*}
    \frac{1}{16-6 + 2\sqrt{5}} = \frac{1+x^2}{x^2 + b^2-1 + 2x \sqrt{b^2-1}}
\end{equation*}
\begin{equation}
    \label{EQ: 5 Lighthouse, x fraction}
    \frac{1}{16-6 + 2\sqrt{5}} = \frac{1+x^2}{x^2 + b^2-1 + 2x \sqrt{40 + \frac{5}{4}(6 - 2\sqrt{5}) -15\sqrt{6 - 2\sqrt{5}} }}
\end{equation}
Again, we were unable to reduce this and used Wolfram$|$Alpha \cite{wolfram-n5-edge} to calculate $x \approx 1.2471$. Measuring $x$ using GeoGebra also gives $x \approx 1.2471$. $D(5)$ is then given as
\begin{equation}
    D(5) = 5(1.2471 - \arctan(1.2471)) \approx 1.7609 
\end{equation}
\subsection{A general rule for any number of lighthouses?}
\begin{figure}[h!]
    \centering
    \begin{minipage}{0.48\textwidth}
        \centering
        \captionsetup{width=.9\textwidth}
        \includegraphics[width=0.9\textwidth, frame]{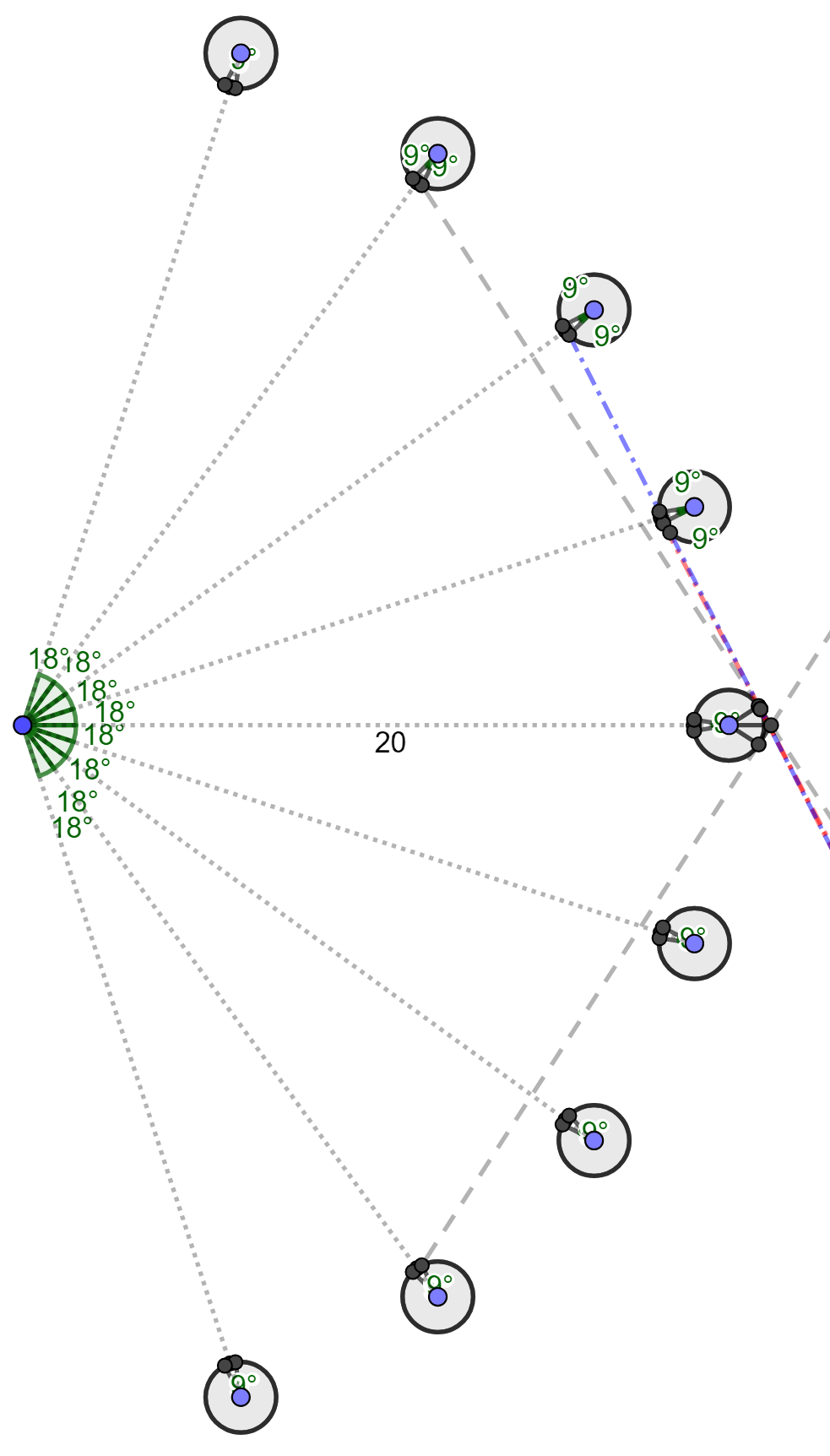} % first figure itself
        \caption{20 Lighthouses, 9 of them are drawn.}
        \label{FIG: 20 Lighthouse - Edge}
    \end{minipage}
    \hfill
    \begin{minipage}{0.48\textwidth}
        \centering
        \captionsetup{width=.9\textwidth}
        \includegraphics[width=0.9\textwidth, frame]{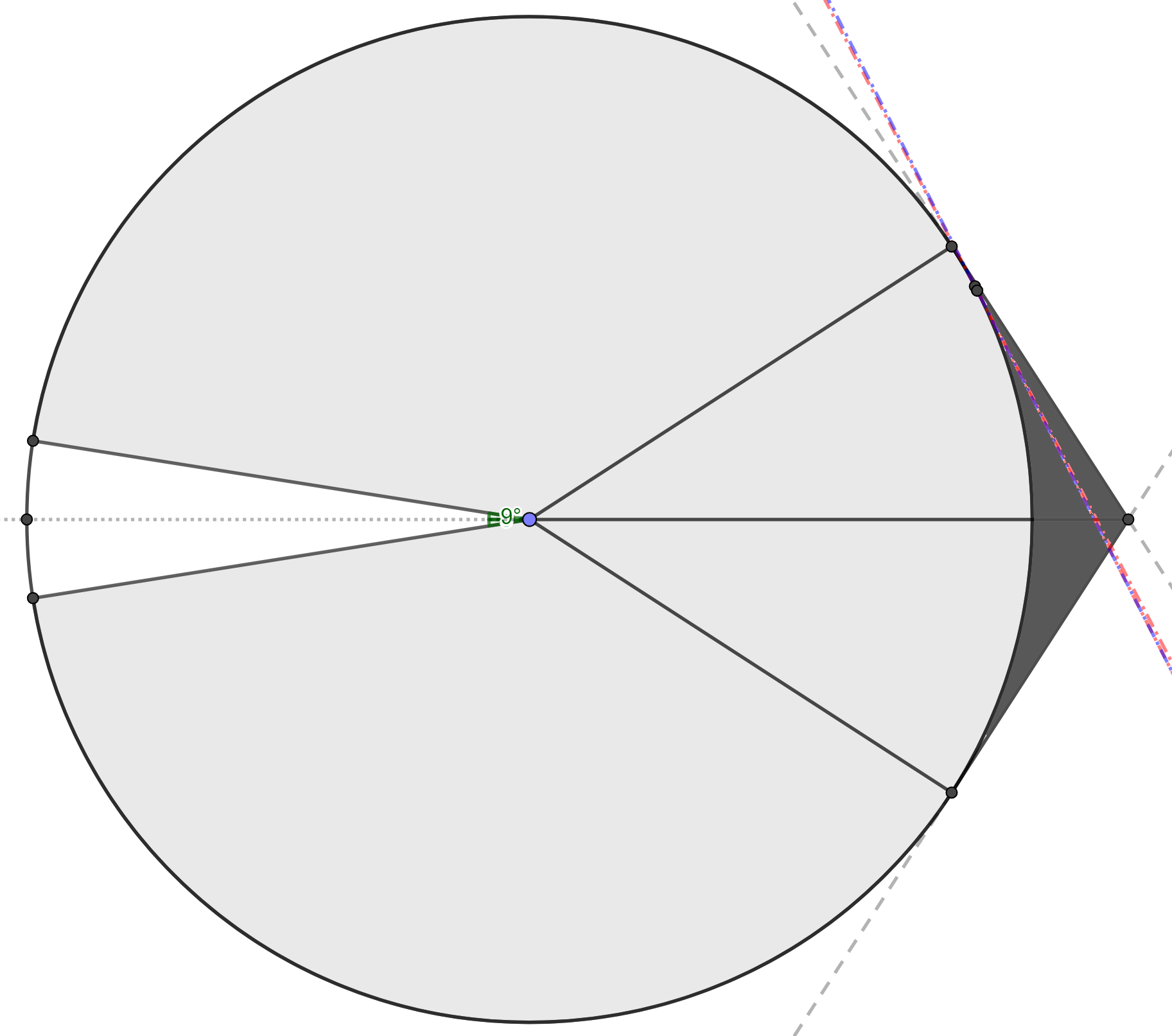} % second figure itself
        \caption{20 Lighthouses, zoomed on target lighthouse. Red striped line is coming from the closest lighthouse, blue striped line is coming from the second closest lighthouse, the gray striped line, which actually defines the dark area, is coming from the third closest lighthouse.}
        \label{FIG: 20 Lighthouse - Edge - Zoom}
    \end{minipage}
\end{figure}

\begin{figure}
    \centering
    \includegraphics[width=0.75\textwidth]{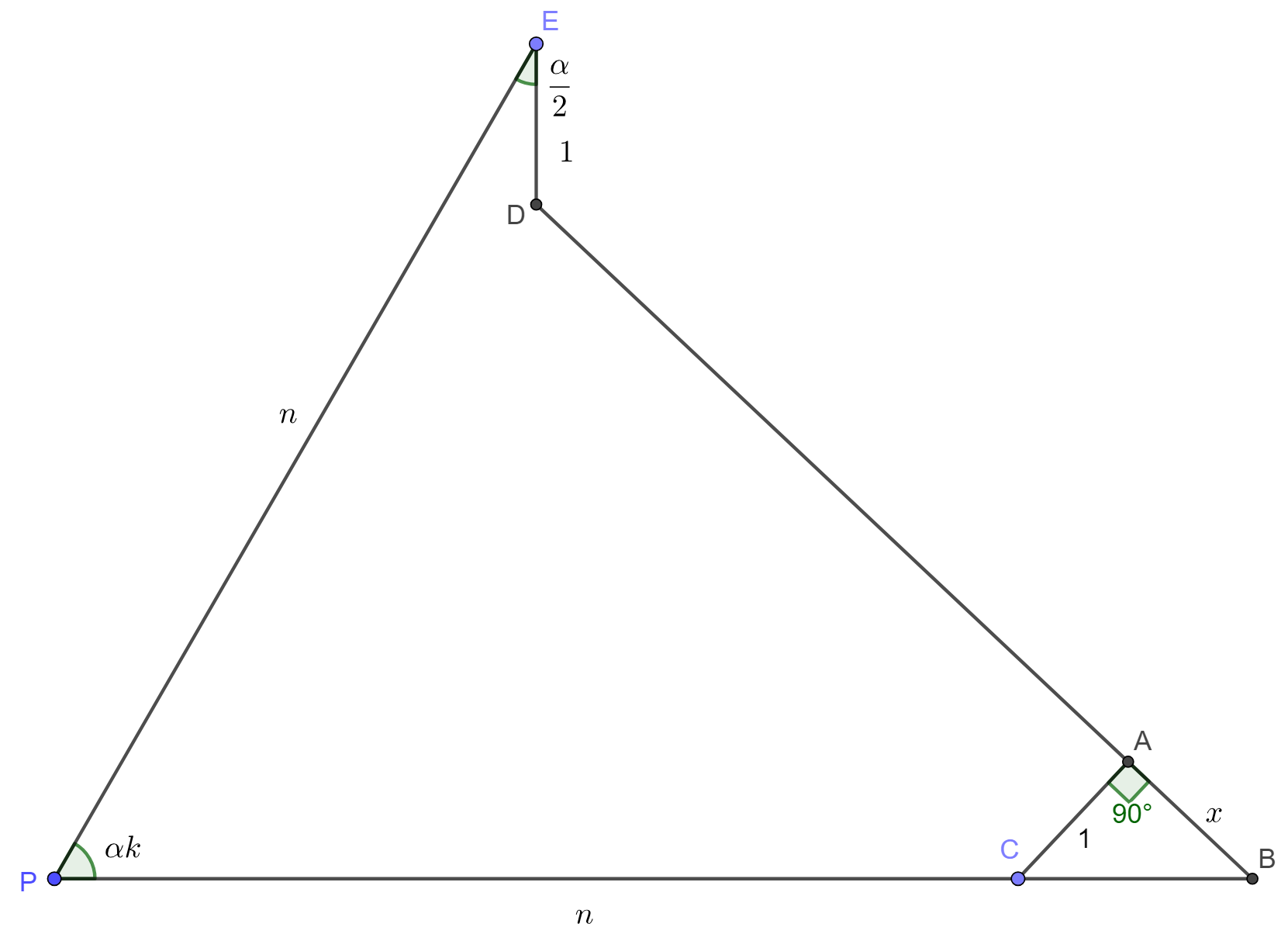} % second figure itself
    \caption{The geometry between a target lighthouse and the source lighthouse that illuminates behind it}
    \label{FIG: N Lighthouse - Edge - General}
\end{figure}

So far, we have failed to reduce the equation obtained from the similarity of the triangles, we sought help from Wolfram$|$Alpha and compared the result to the measurement using GeoGebra. Nevertheless, we could try to find a general form for $x$, basing a claim that the dark area is defined by the closest lighthouses of the target lighthouse, like we did for the previous variation where the dark area was defined by the furthest lighthouses when $n$ was odd. But we realize this is not the case. Upon generalization, we might think that maybe the illumination angle $\alpha$ becomes so small that the closest neighbor is no more able to illuminate behind the target lighthouse, instead, some other lighthouse pair does the job. To approach this case, consider a target lighthouse $L_0$ and a source lighthouse $L_s$. At most, the light ray coming out of the point on the arc of $L_s$ would be tangent to $L_s$. More than that would mean that the light ray actually passes through $L_s$ to be tangent to $L_0$. This is exactly the case for $n=20$. At $n=19$ the angle is slightly more than $90^\circ$ but at $n=20$ the angle is less than $90^\circ$, so perhaps some other lighthouse is illuminating the target lighthouse. For $n=20$, the $3^{rd}$ closest lighthouse is the one illuminating the target lighthouse, as shown in figures 
\ref{FIG: 20 Lighthouse - Edge} and \ref{FIG: 20 Lighthouse - Edge - Zoom}.
So even if we had a formula to calculate $x$ based on the cases $n=3,4,5$ it would cease to work after $n=19$. The general figure is given by \ref{FIG: N Lighthouse - Edge - General}.
Knowing that $\alpha = 360^\circ/n$, describe $x$ and $k$ in terms of $n$, that is the task! $k$ is the number of the lighthouse, $k=1$ means the closest lighthouse, $k=2$ means the second closest lighthouse, $k=3$ the third closest and so on... For example, $k=1$ for values $n=1,2,...,19$ and then $k=3$ for $n=20$. After we actually find a way to deal with $k$ we will find $x$ (which the author had trouble even without $k$). The case of $n=20$ was found empirically, and we do not have a formula that tells us which lighthouse will illuminate the target lighthouse for a given number of lighthouses.
%% END %%

%% file: Sections/3-Conclusion.tex
%% BEGIN %%
\section{Conclusion \& Closing Remarks}
In this section we will recapitulate our findings on both variations of the problem.

If the lighting is done using a single point light source at the center of each lighthouse, we can define a piece-wise formula for the dark area.
\begin{definition}[Total Dark area when there is a point light source at the center]
\begin{equation}
\label{EQ: Conclusion on Center}
D(n) = 
\begin{cases} 
      0 & n = 1 \\
      \infty & n \equiv 0 \pmod{2} \\
      n(x_n - \arctan(x_n)) & n \equiv 1 \pmod{2}, n > 1  
\end{cases}
\end{equation}
\begin{equation}
x_n = \frac{\sqrt{4n^2\cos^2\left(\frac{\pi}{2n}\right) - 1} + 2n^2\sin(\frac{\pi}{n})\cos^2(\frac{\pi}{2n})}{ n^2\sin^2\left(\frac{\pi}{n}\right)-1}
\end{equation}
\end{definition}
We would like to ask, what is $\lim_{n \to \infty}n(x_n - \arctan(x_n))$? Furthermore, as $n$ goes to infinity we have shown that the distance between two lighthouses is $2\pi$, on the other hand $\alpha = 2\pi/n$ so the illumination angle would be approaching 0. Reconsidering figure \ref{FIG: N Odd Lighthouses - Center} we would say that $y=\pi$ but then we have the angle $\angle CPE$ approaching $180^\circ$, so does the lighthouse on the left illuminate the target lighthouse on the right at all? We believe that it does illuminate the target lighthouse but the tangent light rays are almost parallel to each other so the dark area approaches infinity, but we do not have a proof for this yet. Furthermore, there could be a better way to calculate $x_n$ instead of brute-forcing our way to $x$ with Cosine theorems and Pythagoras rules.

As for the point light source at arc case, there is another problem regarding which lighthouse is illuminating the target lighthouse. We have to find a rule regarding which lighthouse illuminates the target lighthouse for given $n$, which enables us to write $k$ in terms of $n$, then we could find $x$ in terms of $n$.

The results obtained in this paper can be seen on table \ref{TABLE: Findings}. The question remains: What is the total dark area for any $n$ for both variations of the problem, and what is it as $n$ approaches infinity?

\begin{table}[h]
\caption{Results obtained in this paper for the total dark area for a given number of lighthouses.}
\begin{center}
\begin{tabular}{|c|c|c|}
\hline
\begin{tabular}[c]{@{}c@{}}Number of\\ lighthouses\end{tabular}  & \begin{tabular}[c]{@{}c@{}}Point Light Source\\ at the Center\end{tabular} &   \begin{tabular}[c]{@{}c@{}}Point Light Sources\\ at the Arc\end{tabular} \\
\hline
$1$  &  $0$ &  $0$       \\ \hline
$2$ &  $\infty$ &   $\infty$    \\ \hline
$3$ & $5.0376$ &    $3.4665$  \\ \hline 
$4$ & $\infty$ &  $2.2475$      \\ \hline
$5$ & $16.7851$ & $1.7609$       \\ \hline
$\vdots$ & $\vdots$ & $\vdots$     \\ \hline
$n$ & \begin{tabular}[c]{@{}c@{}}Equation \eqref{EQ: Conclusion on Center} for odd $n$\\  $\infty$ for even $n$\end{tabular} & ?     \\ \hline
$\infty$ & Possibly $\infty$ & ?     \\
\hline
\end{tabular}
\end{center}
\label{TABLE: Findings}
\end{table}

%% END %%